\documentclass[11pt, a4paper]{article}
\usepackage[letterpaper, left=2cm, right=2cm, top=2cm, bottom=2cm]{geometry}
\usepackage{amsmath}
\usepackage{graphicx}
\usepackage{dsfont}
\usepackage{amsthm}
\usepackage{color}
\usepackage{appendix}
\usepackage{url}
\usepackage{wrapfig}
\usepackage{subfigure}
\usepackage{multirow}
\usepackage{xcolor}
\usepackage{cancel}
\usepackage{psfrag}
\usepackage{caption}
\usepackage{soul}

\usepackage{subeqnarray}
\usepackage{img_macros}

\title{Solution landscapes in nematic microfluidics.}
\author{ M. CRESPO$^{1,\ast}$, I.~M.~GRIFFITHS$^2$ A.~MAJUMDAR$^3$ and A. M RAMOS$^1$ \\[0.3cm]
$^1$ {\small Departamento de Matem\'atica Aplicada, Universidad Complutense de Madrid} \\[0.1cm]
{\small$\&$ Instituto de Matem\'atica Interisciplinar}\\[0.1cm]
{\small Plaza de Ciencias, 3, 28040 Madrid, Spain.}\\[0.1cm]
{\small $^\ast$ E-mail: mcresp01@ucm.es. Tel.: +34-913944462} \\[0.3cm]
$^2$ {\small Mathematical Institute, Radcliffe Observatory Quarter,}\\[0.1cm]
{\small University of Oxford, Oxford, OX2 6GG, U.K.}\\[0.3cm]
$^3$ {\small Department of Mathematical Sciences, University of Bath,}\\[0.1cm]
{\small Bath, BA2 7AY, U.K.}}

\date{18 July 2016}
\begin{document}
\maketitle
\begin{abstract}
We study the static equilibria of a simplified Leslie--Ericksen model for a unidirectional uniaxial nematic flow in a prototype microfluidic channel, as a function of the pressure gradient $\mathcal{G}$ and inverse anchoring strength, $\mathcal{B}$. We numerically find multiple static equilibria for admissible pairs $(\mathcal{G}, \mathcal{B})$ and classify them according to their winding numbers and stability. The case $\mathcal{G}=0$ is analytically tractable and we numerically study how the solution landscape is transformed as $\mathcal{G}$ increases. We study the one-dimensional dynamical model, the sensitivity of the dynamic solutions to initial conditions and the rate of change of $\mathcal{G}$ and $\mathcal{B}$. We provide a physically interesting example of how the time delay between the applications of $\mathcal{G}$ and $\mathcal{B}$ can determine the selection of the final steady state.
\end{abstract}
\section{Introduction}\label{sec:Introduction}
Recent years have seen a tremendous surge in research in complex
fluids, of which nematic liquid crystals (NLC) are a prime
example.\cite{dg,senguptathesis,Sengupta2013} Nematic liquid
crystals are anisotropic liquids that combine the fluidity of
liquids with the orientational order of solids i.e. the
constituent rod-like molecules typically align along certain
preferred or distinguished directions and this orientational
anisotropy can have a profound optical signature.\cite{sluckin2004crystals} Microfluidics
is a thriving field of research; scientists typically manipulate
fluid flow, say conventional isotropic fluids, in narrow channels
complemented by different boundary treatments, leading to novel
transport and mixing phenomena for fluids and potentially new
health and pharmaceutical applications.\cite{microfluidics,stone2004engineering,whitesides2006} A
natural question to ask is what happens when we replace a
conventional isotropic liquid with an anisotropic liquid, such as
a nematic liquid crystal?\cite{Sengupta2013} Nematic microfluidics
have recently generated substantial interest by virtue of their
optical, rheological and backflow properties along with their
defect profiles.\cite{fancyref}

In Sengupta et al.,\cite{Sengupta2013} the authors investigate, both
experimentally and numerically, microfluidic channels filled with
nematic solvents.  The authors work with a thin microfluidic
channel with length much greater than width and width much greater
than depth. A crucial consideration is the choice of boundary
conditions and the authors work with homeotropic or normal
boundary conditions on the top and bottom channel surfaces, which
require the molecules to be oriented in the direction of the
surface normal. The anchoring strength is a measure of how
strongly the boundary conditions are enforced: strong anchoring
roughly corresponds to Dirichlet conditions for the director field
and zero anchoring describes free (Neumann homogeneous) boundary
conditions. We expect most experiments to have moderate to strong
anchoring conditions. The authors impose a flow field transverse
to the anchoring conditions so that there are at least two
competing effects in the experiment: anchoring normal to the
boundaries and flow along the length of the microfluidic channel.
They work with weak, medium, and strong flow speeds in qualitative
terms and observe complex flow transitions. In the weak-flow
regime, the molecules are only weakly affected by the flow and the
molecular orientations are largely determined by the anchoring
conditions. As the flow strength increases, a complex coupling
between the molecular alignments and the flow field emerges and
the nematic molecules reorient to align somewhat with the flow
field. The medium-flow director field exhibits boundary layers
near the centre and the boundaries where the director field is
strongly influenced by either the flow field or the boundary
conditions. In the strong-flow regime, the molecules are almost
entirely oriented with the flow field, with the exception of thin
boundary layers near the channel surfaces to match the boundary
conditions. The authors study these transitions experimentally and
their experimental results suggest a largely \textit{uniaxial}
profile wherein the molecules exhibit a single distinguished
direction of molecular alignment and this direction is referred to
as being the \textit{director} in the literature.\cite{p1995physics} The authors present experimental measurements for the optical profiles and flow fields and their experimental work is
complemented by a numerical analysis of the nematodynamic
equations in the Beris--Edwards theory.\cite{yeomans} The
Beris--Edwards theory is one of the most general formulations of
nematodynamics, that accounts for both uniaxial and biaxial
systems (with a primary and secondary direction of molecular
alignment) and variations in the degree of orientational order. The authors   numerically
reproduce the experimentally observed flow transitions, the
 director and flow-field profiles, all of which are
in good qualitative agreement with the experiments.

In Anderson et al.,\cite{Linda2015} the authors model
this experimental set-up within the Leslie--Ericksen model for
nematodynamics. Their Leslie--Ericksen model is restricted to
uniaxial nematics with constant ordering (a constant degree of
orientational order).\cite{ericksen} They present governing
equations for the flow field and the nematic director field along
with the constitutive relations that describe the coupling between
the director and the flow field (see Appendix A for details) and assume that all dependent variables only vary along the channel depth, with a unidirectional flow  along the channel length, consistent with
the experiments. These assumptions greatly simplify the
mathematical model, yielding a decoupled system of partial
differential equations for the director field, which captures the
flow dynamics through a single variable: the pressure gradient, 
$\mathcal{G}$, along the channel length. The authors define two
separate boundary-value problems: one for weak-flow solutions and
one for strong-flow solutions, described by two different sets of
boundary conditions for the director field. They find weak- and
strong-flow solutions for all values of the pressure gradient and
they relate the resulting flow profile to the mean flow speed by a
standard Poiseuille-flow-type relation. 
The energy of the weak-flow solution is lower than the strong-flow solution for small $\mathcal{G}$ and
there is an energy cross-over at some critical value,
$\mathcal{G}^{\ast}$, that depends on the anchoring strength at the
channel surfaces.

In this paper, we build on the work in Anderson et al.\cite{Linda2015} by performing an extensive study of the static solution landscape, complemented by some
numerical investigations of the dynamical behavior, as the system
evolves to these equilibrium configurations. We adopt the same
model with the same underpinning assumptions as in Anderson et al.,\cite{Linda2015} but we do not define two separate boundary-value
problems. We impose weak anchoring conditions for the director
field on the top and the bottom surfaces since it
includes both the weak and strong anchoring configurations and allow us to capture the competition between the flow field and the anchoring strength. 

We compute the static equilibrium solutions, using a combination of
analytic and numerical methods, as a function of $\mathcal {G}$
and the inverse anchoring strength $\mathcal{B}$. The case $\mathcal{G}=0$ is analytically tractable and we identify two different classes of solutions and characterize their stability. This is complemented by
an asymptotic analysis in the limits $\mathcal{G}\to 0$ and
$\mathcal{G}\to \infty$, with the latter regime yielding useful
information about the boundary layers near channel surfaces, which are
experimentally observed in the strong-flow regimes.\cite{Sengupta2013} We then study the solution landscape for $\mathcal{G}\neq 0$ and track the stable and unstable solution branches as a function of $\left(\mathcal{G}, \mathcal{B} \right)$. Our work largely focuses on the static
equilibria but the last section is devoted to a numerical study
of the dynamic Leslie--Ericksen model and its sensitivity to the initial condition.
In particular, we present a numerical example for which we can control the final steady state by manipulating the rate of change of the
pressure gradient and anchoring conditions.

The paper is organized as follows. In
Section~\ref{sec:mathematicalmodel}, we present the
Leslie--Ericksen dynamic model, the governing equations and
boundary conditions. In Section~\ref{sec:equilibrium}, we explore
the static solution landscape as a function of the pressure
gradient and anchoring strength. In Section~\ref{sec:time}, we
study the dynamic model, with focus on the effects of initial
conditions and the time-dependent forms of the pressure gradient
and anchoring strength, and conclude in
Section~\ref{sec:conclusions} by putting our work in context and
discuss future developments. 

\section{Mathematical Model}\label{sec:mathematicalmodel}
\begin{figure}[ht!]
\centering
%
%
\begin{psfrags}%
\psfragscanon%
%
\psfrag{s05}[l][l]{\color[rgb]{0,0,0}\setlength{\tabcolsep}{0pt}\begin{tabular}{l}{\Huge $\hat{z}$}\end{tabular}}%
\psfrag{s06}[l][l]{\color[rgb]{0,0,0}\setlength{\tabcolsep}{0pt}\begin{tabular}{l}{\Huge $\hat{x}$}\end{tabular}}%
\psfrag{s07}[l][l]{\color[rgb]{0,0,0}\setlength{\tabcolsep}{0pt}\begin{tabular}{l}{\Huge $-h$}\end{tabular}}%
\psfrag{s08}[l][l]{\color[rgb]{0,0,0}\setlength{\tabcolsep}{0pt}\begin{tabular}{l}{\Huge $h$}\end{tabular}}%
\psfrag{s09}[l][l]{\color[rgb]{1,0,0}\setlength{\tabcolsep}{0pt}\begin{tabular}{l}{\Huge Fluid Flow}\end{tabular}}%
%
\psfrag{x01}[t][t]{0}%
\psfrag{x02}[t][t]{0.1}%
\psfrag{x03}[t][t]{0.2}%
\psfrag{x04}[t][t]{0.3}%
\psfrag{x05}[t][t]{0.4}%
\psfrag{x06}[t][t]{0.5}%
\psfrag{x07}[t][t]{0.6}%
\psfrag{x08}[t][t]{0.7}%
\psfrag{x09}[t][t]{0.8}%
\psfrag{x10}[t][t]{0.9}%
\psfrag{x11}[t][t]{1}%
\psfrag{x12}[t][t]{0}%
\psfrag{x13}[t][t]{0.1}%
\psfrag{x14}[t][t]{0.2}%
\psfrag{x15}[t][t]{0.3}%
\psfrag{x16}[t][t]{0.4}%
\psfrag{x17}[t][t]{0.5}%
\psfrag{x18}[t][t]{0.6}%
\psfrag{x19}[t][t]{0.7}%
\psfrag{x20}[t][t]{0.8}%
\psfrag{x21}[t][t]{0.9}%
\psfrag{x22}[t][t]{1}%
%
\psfrag{v01}[r][r]{0}%
\psfrag{v02}[r][r]{0.1}%
\psfrag{v03}[r][r]{0.2}%
\psfrag{v04}[r][r]{0.3}%
\psfrag{v05}[r][r]{0.4}%
\psfrag{v06}[r][r]{0.5}%
\psfrag{v07}[r][r]{0.6}%
\psfrag{v08}[r][r]{0.7}%
\psfrag{v09}[r][r]{0.8}%
\psfrag{v10}[r][r]{0.9}%
\psfrag{v11}[r][r]{1}%
\psfrag{v12}[r][r]{0}%
\psfrag{v13}[r][r]{0.1}%
\psfrag{v14}[r][r]{0.2}%
\psfrag{v15}[r][r]{0.3}%
\psfrag{v16}[r][r]{0.4}%
\psfrag{v17}[r][r]{0.5}%
\psfrag{v18}[r][r]{0.6}%
\psfrag{v19}[r][r]{0.7}%
\psfrag{v20}[r][r]{0.8}%
\psfrag{v21}[r][r]{0.9}%
\psfrag{v22}[r][r]{1}%
%
\resizebox{5cm}{!}{\includegraphics{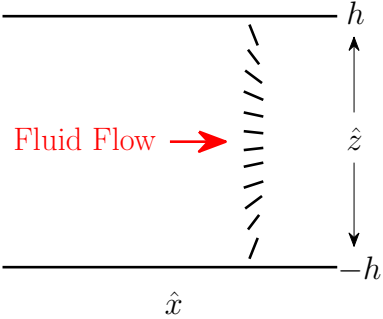}}%
\end{psfrags}%
%

\caption{Schematic of the microfluidic channel set-up. The nematic molecules are anchored at the top and bottom surfaces and are deformed by the fluid flow from the left.}
\label{fig:scheme}
\end{figure}

As in Anderson et al.,\cite{Linda2015} we model the NLC within the microfluidic channel in the Leslie--Ericksen framework. The channel has dimensions, $L_{\hat{x}} >> L_{\hat{y}} >> L_{\hat{z}}$, in the $\hat{x}$, $\hat{y}$ and $\hat{z}$ directions respectively, consistent with the experimental set-up in Anderson et al.\cite{Linda2015} and Sengupta et al.\cite{Sengupta2013} The NLC is purely uniaxial with constant order parameter, by
assumption, and is hence fully described by a director field,
$\ve{n}$, that represents the single preferred direction of
nematic alignment. Here, $\ve{n}$ and $-\ve{n}$ are physically indistinguishable (in the absence of polarity the sign of $\ve{n}$ has no physical meaning). 
We additionally assume that all dependent variables only depend
on the $\hat{z}$-coordinate, along the channel depth, as depicted
in Figure \ref{fig:scheme}.
Then the director field is of the form
$\ve{n}=(\sin(\theta(\hat{z},\hat{t})),0,\cos(\theta(\hat{z},\hat{t})))$
and  the velocity field is unidirectional, of the form
$\ve{v}=(u(\hat{z},\hat{t}),0,0)$, with $-h\leq \hat{z} \leq h$.
Since $\ve{n}$ and $-\ve{n}$ are indistinguishable, $\theta$ and $\theta+k\pi$, $k\in \mathds{Z}$, describe the same director profile. We assume that $u(\hat{z},\hat{t})$ is symmetric around the
center-line (i.e around $\hat{z}=0$) and no-slip conditions are imposed
on the channel walls (i.e. $u(\pm h,\hat{t})=0$). We assume weak anchoring boundary conditions for $\theta$ on
$\hat{z}= \pm h$, that can be derived from the well-known Rapini--Papoular weak-anchoring energy,\cite{Rapini1969}
 $$ E_S = \int_{\hat{z}=\pm h} \frac{A}{2}\sin^2 \theta~ d\hat{x}~d\hat{y}, $$
which enforces $\theta(-h)=k_{1}\pi$ and $\theta(h)=k_{2}\pi$ ($k_{1},k_{2}\in \mathds{Z}$) for large anchoring coefficients $A>0$. In other words, the  Rapini--Papoular energy enforces homeotropic anchoring (along the normal to the surface) described by, $\ve{n}= \pm \left(0, 0, 1 \right)$ on $\hat{z} = \pm h$.
 
 We substitute the assumed forms for $\ve{v}$ and $\ve{n}$
into the full governing equations, as outlined in \ref{ericksenmodel}, and using (\ref{ledynamic}\textit{a,b}), we obtain the following decoupled initial-boundary-value problem for $\theta$:
\begin{align}\label{ledynamictheta}
\left\{
\begin{array}{l c l}
\big(\hat{\gamma}_{1}\hat{g}(\theta) - \hat{m}(\theta)^{2}\big) \displaystyle \frac{\partial \theta}{\partial \hat{t}} =  K\hat{g}(\theta)\displaystyle\frac{\partial^{2} \theta}{\partial \hat{z}^{2}} +G\hat{z}\hat{m}(\theta) & & \hat{z}\in (-h,h),\,\, \hat{t}>0,\\
\\
\theta(\hat{z},0) =\Theta(\hat{z}) & & \hat{z} \in (-h,h),
\\
\\
K \displaystyle \frac{\partial \theta}{\partial \hat{z}}= - \frac{A}{2}\sin(2\theta(\hat{z},\hat{t})) & & \hat{z} =h,\,\, \hat{t}>0,
\\
\\
K \displaystyle \frac{\partial \theta}{\partial \hat{z}}= \frac{A}{2}\sin(2\theta(\hat{z},\hat{t})) &  & \hat{z} =-h,\,\, \hat{t}>0,
\end{array}\right.
\end{align}
where subscripts denote partial differentiation, $K$ (N) is the
elastic constant of the NLC, $\Theta$ is
the initial condition, $-G= \frac{\partial P}{\partial \hat{x}}$ is the component of the
pressure gradient in the channel direction  and $A$ (Nm$^{-1}$) is the
surface anchoring strength. Note that for a physically realistic solution, we expect that as $A \to \infty$, $2\theta$ tends to an integer multiple of $\pi$ on $\hat{z}=\pm h$.
The functions
$$\hat{m}(\theta)=  \hat{\alpha}_{2}\cos^{2}(\theta)-\hat{\alpha}_{3}\sin^{2}(\theta) \mbox{ and}$$
$$\hat{g}(\theta) = \hat{\alpha}_{1}\cos^{2}(\theta)\sin^{2}(\theta)  +
\frac{\hat{\alpha}_{5}-\hat{\alpha}_{2}}{2}\cos^{2}(\theta) +
\frac{\hat{\alpha}_{3}+\hat{\alpha}_{6}}{2} \sin^{2}(\theta)+
\frac{\hat{\alpha}_{4}}{2},$$ the $\hat{\alpha}_{i}$ (N m$^{-2}$
s), $i\in \{1, \ldots,6\}$, are constant viscosities related to
each other by the Parodi relation,\cite{Parodi}
$\hat{\alpha}_{2}+\hat{\alpha}_{3}=\hat{\alpha}_{6}-\hat{\alpha}_{5}$.
Characteristic values for the dimensionless nematic viscosities
are $\alpha_{1}=-0.1549$, $\alpha_{\rm 2}=-0.9859$,
$\alpha_{3}=-0.0535$, $\alpha_{5}=0.7324$ and $\alpha_{\rm
6}=-0.39$.\cite{Linda2015} Note that the following inequalities must be satisfied
(see Appendix \ref{rmk:coef}):
\begin{equation}\label{strictinequalities}
\hat{g}(\theta)>0, \hspace{2cm} \hat{\gamma}_{1}\hat{g}(\theta)>\hat{m}^{2}(\theta),
\end{equation}
where $\hat{\gamma}_{1}=\hat{\alpha}_{3}-\hat{\alpha}_{2}$.\\
Note that, if $\theta_{1}$, $\theta_{2}$ are the solutions of \eqref{ledynamictheta} corresponding, respectively, to initial conditions $\Theta_{1}(\hat{z})$ and $\Theta_{2}(\hat{z})=\Theta_{1}(\hat{z})+k\pi$ ($k\in \mathds{Z}$), then $\theta_{2}= \theta_{1}+k\pi$ and both $\theta_{1}$ and $\theta_{2}$ correspond  to the same physical description of molecular orientation.\\

We non-dimensionalize the system \eqref{ledynamictheta} using the scalings
$$z= \frac{\hat{z}}{h}, \hspace{1cm} \alpha_{i}= \frac{\hat{\alpha}_{i}}{\hat{\alpha}_{4}}, \hspace{1cm}  \gamma_{1}= \frac{\hat{\gamma}_{1}}{\hat{\alpha}_{4}}, \hspace{1cm} t= \frac{K\hat{t}}{\hat{\alpha}_4h^2}. $$
The dimensionless version of
(\ref{ledynamictheta}) is then
\begin{subequations}\label{systemtdep}
\begin{align}
\big(\gamma_{1}g(\theta)-m(\theta)^2\big)\displaystyle \frac{\partial \theta}{\partial t} = g(\theta)\displaystyle \frac{\partial^{2} \theta}{\partial z^{2}} + \mathcal{G} z m(\theta) &  & z\in (-1,1),\,\, t>0,\\
\theta(z,0) =\Theta(z) & & z \in (-1,1),\\
\mathcal{B}\displaystyle \frac{\partial \theta}{\partial z}(1,t)=-\sin(2\theta(1,t)) &  & t>0, \label{boundtop}\\
\mathcal{B}\displaystyle \frac{\partial \theta}{\partial z}(-1,t)=\sin(2\theta(-1,t)) & & t>0, \label{boundbot}
\end{align}
\end{subequations}
where $\mathcal{G}={h^{3}G}/{K}$ and $\mathcal{B}={2K}/{Ah}$ are
the dimensionless pressure gradient and the dimensionless inverse
anchoring strength respectively,
$$m(\theta)= \alpha_{2}\cos^{2}(\theta)-\alpha_{3}\sin^{2}(\theta) \mbox{ and }$$
$$g(\theta)  = \alpha_{1}\cos^{2}(\theta)\sin^{2}(\theta)  +\frac{1}{2} \big( (\alpha_{5}-\alpha_{2})\cos^{2}(\theta) + (\alpha_{3}+\alpha_{6})\sin^{2}(\theta)+ 1\big).$$
We compute equilibrium solutions and dynamic time-dependent
solutions of system (\ref{systemtdep}) for different values of
dimensionless pressure gradient $\mathcal{G}$, dimensionless
inverse anchoring strength $\mathcal{B}$ and initial conditions
$\Theta$, using parameter values for the NLC 5CB as
in Anderson et al.\cite{Linda2015}   


\section{Equilibrium Solutions}\label{sec:equilibrium}
We begin by studying the static equilibria of the system (\ref{systemtdep}), $\theta^{\ast}(z)$, which satisfy
\begin{equation}\label{systemstatic}
\left \{\begin{array}{l r}
g(\theta^{\ast}(z))\displaystyle\frac{{\rm d}^{2} \theta^{\ast}}{{\rm d} z^{2}}(z) = -\mathcal{G} z m(\theta^{\ast}(z)) & z\in (-1,1),\\
\\
\mathcal{B}\displaystyle \frac{{\rm d}\theta^{\ast}}{{\rm d}z }(1)=-\sin(2\theta^{\ast}(1)), & \\
\\
\mathcal{B}\displaystyle \frac{{\rm d}\theta^{\ast}}{{\rm d}z }(-1)=\sin(2\theta^{\ast}(-1)). &
\end{array}\right.
\end{equation}
We characterize the equilibrium solutions in terms of their winding number, defined to be
\begin{equation}\label{defn:winding}
\omega(\theta^{\ast})= \frac{\theta^{\ast}(1)-\theta^{\ast}(-1)}{2\pi}.
\end{equation}
The winding number\cite{mermin1979} is a measure of the rotation of the director field between the top and
bottom surfaces. The limit $\mathcal{B}\rightarrow 0$ is the strong anchoring limit, when the boundary conditions on $z=\pm 1$ are strongly enforced and both $\theta^{\ast}(1)$ and $\theta^{\ast}(-1)$ are integer multiples of $\frac{\pi}{2}$ at this limit. Particularly, as we will see in Section \ref{sec:G0}, as $\mathcal{B}\rightarrow 0$, the stable equilibria at $z=\pm 1$ tend to $\theta^{\ast}(\pm 1)= n \pi, n \in \mathds{Z}$ (homeotropic anchoring) and the unstable equilibria to $\theta^{\ast}(\pm 1)= (n +\frac{1}{2})\pi,n \in \mathds{Z}$ (planar anchoring at the boundaries). This is simply because $\theta^{\ast}(\pm 1)= n\pi$ is a minimum of the surface energy used to derive the anchoring conditions at $z=\pm 1$. See \ref{sec:config} for a detailed description of different molecular configurations. In what follows, we track the stable and unstable solutions of \eqref{systemstatic} as the model parameters are varied.
\subsection{No fluid flow ($\mathcal{G}=0$)}\label{sec:G0}
When $\mathcal{G}=0$, we can explicitly solve the system (\ref{systemstatic}) to obtain the static equilibria (see \ref{apsec:steadyg0} for more details). We divide the potentially stable equilibria (see Section \ref{apsec:linearstability}) into two families:
\begin{align}
\mbox{ Type I  		} & \theta_{a_{n}}^{\ast}(z) = a_{n}z, & \mbox{where }  \mathcal{B}a_{n}= -\sin(2a_{n}),  \label{eq:type1m0}
\\ 
\mbox{ Type II  	}&  \theta_{\tilde{a}_{n}}^{\ast}(z) = \tilde{a}_{n}z + \frac{\pi}{2}, & \mbox{where }  \mathcal{B}\tilde{a}_{n}= \sin(2\tilde{a}_{n}). \label{eq:type2m0}
\end{align}
For every value of $\mathcal{B}$, we obtain an ordered set of solutions for \eqref{eq:type1m0}, with $0=a_{0}<a_{1}<\ldots<a_{n}$, $n \in \mathds{N}\cup \{0\}$ depending on $\mathcal{B}$. Moreover, if $a_{n}$ defines a solution, so does $-a_{n}$, which we denote by $a_{-n}$ (identical remarks apply to (\ref{eq:type2m0})). Let $\theta^{\ast}_{a_{n}}$ denote the solution corresponding to $a_{n}$ in \eqref{eq:type1m0}, then $\theta^{\ast}_{a_{n}}= - \theta^{\ast}_{a_{-n}}$ and $\omega(\theta^{\ast}_{a_{n}})= - \omega(\theta^{\ast}_{a_{-n}})= \frac{a_{n}}{\pi}$, where $\omega(\theta^{\ast}_{a_{n}})$ satisfies the transcendental equation 
\begin{equation}\label{eq:betatype1}
\mathcal{B} = - \frac{\sin(2\pi\omega(\theta^{\ast}_{a_{n}}))}{\pi \omega(\theta^{\ast}_{a_{n}})}.
\end{equation}  
Analogous statements apply to solutions $\theta^{\ast}_{\tilde{a}_{n}}$ with $\tilde{a}_{n}$ a solution of equation \eqref{eq:type2m0}, where $\omega(\theta^{\ast}_{\tilde{a}_{n}})$ satisfies the transcendental equation 
\begin{equation}\label{eq:betatype2}
\mathcal{B} =  \frac{\sin(2\pi\omega(\theta^{\ast}_{\tilde{a}_{n}}))}{\pi \omega(\theta^{\ast}_{\tilde{a}_{n}})}.
\end{equation} 
Thus there is a symmetric (with respect to $\omega(\theta^{\ast})=0$) arrangement of solutions, which is physically reasonable since we do not expect to have a preferred twist direction
when $\mathcal{G} = 0.$ In Section \ref{apsec:linearstability} we analyze the linear stability of the equilibria \eqref{eq:type1m0}--\eqref{eq:type2m0} to conclude that
$$\begin{array}{c c}
\mbox{Type I 		} & \mbox{ is stable if }n \mbox{ is even and is unstable if }n \mbox{ is odd},\\
\mbox{Type II 		} & \mbox{ is stable if }n \mbox{ is odd and is unstable if }n \mbox{ is even}.\\
\end{array}$$
\begin{figure}[ht!]
\centering
%
%
\begin{psfrags}%
\psfragscanon%
%
\psfrag{s03}[t][t]{\color[rgb]{0,0,0}\setlength{\tabcolsep}{0pt}\begin{tabular}{c}{\LARGE $\omega(\theta^{\ast}_{a_{n}})$}\end{tabular}}%
\psfrag{s04}[b][b]{\color[rgb]{0,0,0}\setlength{\tabcolsep}{0pt}\begin{tabular}{c}{\LARGE $\mathcal{B}$}\end{tabular}}%
\psfrag{s05}[l][l]{\color[rgb]{0,0,0}\setlength{\tabcolsep}{0pt}\begin{tabular}{l}{\Large $\theta_{a_{-4}}^{\ast}$}\end{tabular}}%
\psfrag{s06}[l][l]{\color[rgb]{0,0,0}\setlength{\tabcolsep}{0pt}\begin{tabular}{l}{\Large$\theta_{a_{-3}}^\ast$}\end{tabular}}%
\psfrag{s07}[l][l]{\color[rgb]{0,0,0}\setlength{\tabcolsep}{0pt}\begin{tabular}{l}{\Large $\theta_{a_{-2}}^\ast$}\end{tabular}}%
\psfrag{s08}[l][l]{\color[rgb]{0,0,0}\setlength{\tabcolsep}{0pt}\begin{tabular}{l}{\Large $\theta_{a_{-1}}^\ast$}\end{tabular}}%
\psfrag{s09}[l][l]{\color[rgb]{0,0,0}\setlength{\tabcolsep}{0pt}\begin{tabular}{l}{\Large $\theta_{a_{0}}^\ast$}\end{tabular}}%
\psfrag{s10}[l][l]{\color[rgb]{0,0,0}\setlength{\tabcolsep}{0pt}\begin{tabular}{l}{\Large $\theta_{a_{1}}^\ast$}\end{tabular}}%
\psfrag{s11}[l][l]{\color[rgb]{0,0,0}\setlength{\tabcolsep}{0pt}\begin{tabular}{l}{\Large $\theta_{a_{2}}^\ast$}\end{tabular}}%
\psfrag{s12}[l][l]{\color[rgb]{0,0,0}\setlength{\tabcolsep}{0pt}\begin{tabular}{l}{\Large $\theta_{a_{3}}^\ast$}\end{tabular}}%
\psfrag{s13}[l][l]{\color[rgb]{0,0,0}\setlength{\tabcolsep}{0pt}\begin{tabular}{l}{\Large $\theta_{a_{4}}^\ast$}\end{tabular}}%
\psfrag{s14}[l][l]{\color[rgb]{0,0,0}\setlength{\tabcolsep}{0pt}\begin{tabular}{l}{\LARGE Type I}\end{tabular}}%
%
\psfrag{x01}[t][t]{-2}%
\psfrag{x02}[t][t]{-$\frac{3}{2}$}%
\psfrag{x03}[t][t]{-1}%
\psfrag{x04}[t][t]{-$\frac{1}{2}$}%
\psfrag{x05}[t][t]{0}%
\psfrag{x06}[t][t]{$\frac{1}{2}$}%
\psfrag{x07}[t][t]{1}%
\psfrag{x08}[t][t]{$\frac{3}{2}$}%
\psfrag{x09}[t][t]{2}%
%
\psfrag{v01}[r][r]{0.1}%
\psfrag{v02}[r][r]{{\Large \textcolor{red}{$\mathcal{B}_{4}^{\ast}= \mathcal{B}_{-4}^{\ast}$}}}%
\psfrag{v03}[r][r]{0.3}%
\psfrag{v04}[r][r]{{\Large \textcolor{red}{$\mathcal{B}_{2}^{\ast}= \mathcal{B}_{-2}^{\ast}$}}}%
\psfrag{v05}[r][r]{0.6}%
\psfrag{v06}[r][r]{0.7}%
\psfrag{v07}[r][r]{0.8}%
\psfrag{v08}[r][r]{0.9}%
\psfrag{v09}[r][r]{1}%
%
\resizebox{8cm}{!}{\includegraphics{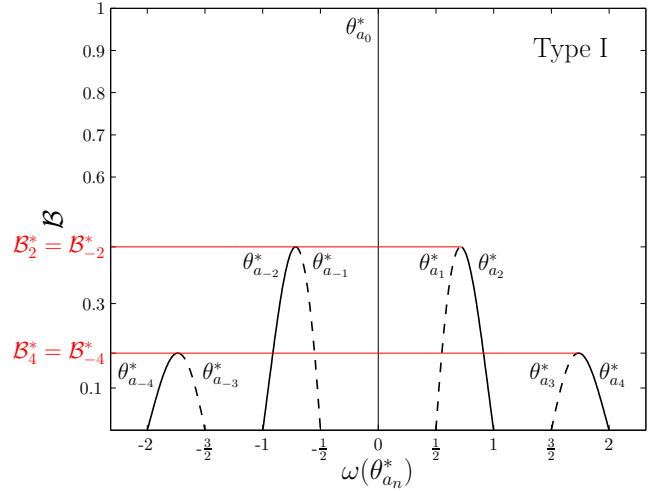}}%
\end{psfrags}%
%

\caption{Case $\mathcal{G}=0$: Solutions of \eqref{eq:betatype1} indicating the emergence of non-constant steady-state solutions $\theta^{\ast}_{a_{n}}$, $n=0,\pm 1, \ldots$ at critical values $\mathcal{B}_{2n}^{\ast}$ for $n=\pm 1, \pm 2,\ldots$. The solid and dashed lines represent, respectively, the values of $\omega(\theta_{a_{n}}^{\ast})$ for which the steady state $\theta_{a_{n}}^{\ast}$ is stable or unstable.}
\label{fig:stabilityg0typ1}
\end{figure}
\begin{figure}
\centering
%
%
\begin{psfrags}%
\psfragscanon%
%
\psfrag{s03}[t][t]{\color[rgb]{0,0,0}\setlength{\tabcolsep}{0pt}\begin{tabular}{c}{\LARGE $\omega(\theta^{\ast}_{\tilde{a}_{n}})$}\end{tabular}}%
\psfrag{s04}[b][b]{\color[rgb]{0,0,0}\setlength{\tabcolsep}{0pt}\begin{tabular}{c}{\LARGE $\mathcal{B}$}\end{tabular}}%
\psfrag{s05}[l][l]{\color[rgb]{0,0,0}\setlength{\tabcolsep}{0pt}\begin{tabular}{l}{\Large $\theta_{\tilde{a}_{-3}}^\ast$}\end{tabular}}%
\psfrag{s06}[l][l]{\color[rgb]{0,0,0}\setlength{\tabcolsep}{0pt}\begin{tabular}{l}{\Large $\theta_{\tilde{a}_{-2}}^\ast$}\end{tabular}}%
\psfrag{s07}[l][l]{\color[rgb]{0,0,0}\setlength{\tabcolsep}{0pt}\begin{tabular}{l}{\Large $\theta_{\tilde{a}_{-1}}^\ast$}\end{tabular}}%
\psfrag{s08}[l][l]{\color[rgb]{0,0,0}\setlength{\tabcolsep}{0pt}\begin{tabular}{l}{\Large $\theta_{\tilde{a}_{0}}^\ast$}\end{tabular}}%
\psfrag{s09}[l][l]{\color[rgb]{0,0,0}\setlength{\tabcolsep}{0pt}\begin{tabular}{l}{\Large $\theta_{\tilde{a}_{1}}^\ast$}\end{tabular}}%
\psfrag{s10}[l][l]{\color[rgb]{0,0,0}\setlength{\tabcolsep}{0pt}\begin{tabular}{l}{\Large $\theta_{\tilde{a}_{2}}^\ast$}\end{tabular}}%
\psfrag{s11}[l][l]{\color[rgb]{0,0,0}\setlength{\tabcolsep}{0pt}\begin{tabular}{l}{\Large $\theta_{\tilde{a}_{3}}^\ast$}\end{tabular}}%
\psfrag{s12}[l][l]{\color[rgb]{0,0,0}\setlength{\tabcolsep}{0pt}\begin{tabular}{l}{\LARGE Type II}\end{tabular}}%
%
\psfrag{x01}[t][t]{-2}%
\psfrag{x02}[t][t]{-$\frac{3}{2}$}%
\psfrag{x03}[t][t]{-1}%
\psfrag{x04}[t][t]{-$\frac{1}{2}$}%
\psfrag{x05}[t][t]{0}%
\psfrag{x06}[t][t]{$\frac{1}{2}$}%
\psfrag{x07}[t][t]{1}%
\psfrag{x08}[t][t]{$\frac{3}{2}$}%
\psfrag{x09}[t][t]{2}%
%
\psfrag{v01}[r][r]{{\Large \textcolor{red}{$\mathcal{B}^{\ast}_{3}=\mathcal{B}^{\ast}_{-3}$} }}%
\psfrag{v02}[r][r]{0.7}%
\psfrag{v03}[r][r]{1}%
\psfrag{v04}[r][r]{1.5}%
\psfrag{v05}[r][r]{{\Large \textcolor{red}{$\mathcal{B}^{\ast}_{1}=\mathcal{B}^{\ast}_{-1}$} }}%
\psfrag{v06}[r][r]{2.3}%
%
\resizebox{8cm}{!}{\includegraphics{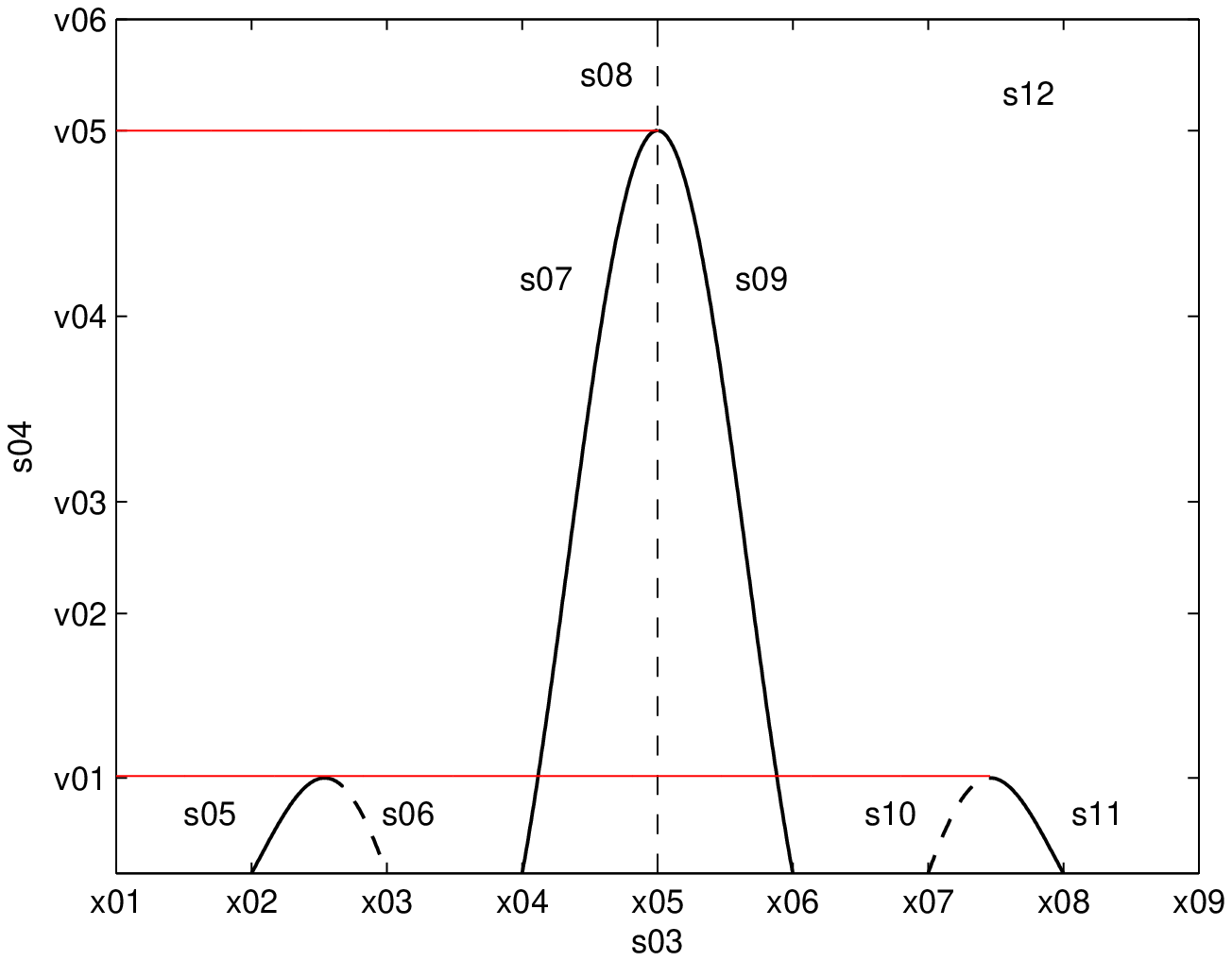}}%
\end{psfrags}%
%

\caption{Case $\mathcal{G}=0$: Solutions of \eqref{eq:betatype2} indicating the emergence of non-constant steady-state solutions $\theta^{\ast}_{\tilde{a}_{n}}$, $n=0,\pm 1, \ldots$ at critical values $\mathcal{B}_{2n+1}^{\ast}$ for $n=0,\pm 1, \ldots$. The solid and dashed lines represent, respectively, the values of $\omega(\theta_{\tilde{a}_{n}}^{\ast})$ for which the steady state $\theta_{\tilde{a}_{n}}^{\ast}$ is stable or unstable.}
\label{fig:stabilityg0typ2}
\end{figure} 
It is clear that the director profiles for  $\theta^{\ast}_{a_{n}}$ and $ \theta^{\ast}_{a_{-n}}$ are reflections of each other about the angle $\theta=0$.
The constant solutions $\theta_{a_{0}}^{\ast} \equiv 0$ and $\theta_{\tilde{a}_{0}}^{\ast} \equiv \frac{\pi}{2}$ exist for all values of $\mathcal{B}$. These are the only solutions for large values of $\mathcal{B}$. Non--constant solutions subject to \eqref{eq:type1m0} and \eqref{eq:type2m0} emerge as $\mathcal{B}$ decreases. \\
We define critical values $\mathcal{B}^{\ast}_{2n}$ with $n = \pm 1, \pm 2, \ldots$ such that, for $n > 0$, the solution branches, $\big(\omega(\theta^{\ast}_{a_{2n}}),\mathcal{B}\big)$ and $\big(\omega(\theta^{\ast}_{a_{2n-1}}),\mathcal{B}\big)$ (and $\big(\omega(\theta^{\ast}_{a_{2n+1}}),\mathcal{B}\big)$ if $n<0$) coalesce at the critical value $\mathcal{B} = \mathcal{B}^{\ast}_{2n}$ and cease to exist for $\mathcal{B} > \mathcal{B}^{\ast}_{2n}$ (see Figure \ref{fig:stabilityg0typ1}). Similarly, we define the critical values  $\mathcal{B}^{\ast}_{2n+1}$ with $n=0, \pm1 , \ldots$ as the coalescence points for solutions of Type II (see Figure \ref{fig:stabilityg0typ2} for a complete description). Solutions with large winding numbers are only observable in the strong--anchoring limit. Notice that for $\mathcal{B}\rightarrow 0$ the stable equilibria are either $\theta_{a_{n}}^{\ast}$ with $\omega(\theta^{\ast}_{a_{n}})=k\pi$ or $\theta_{\tilde{a}_{n}}^{\ast}$ with $\omega(\theta^{\ast}_{\tilde{a}_{n}})=(k+\frac{1}{2})\pi$, $k\in \mathds{Z}$, and in both cases $\theta^{\ast}(\pm 1)$ tends to a multiple of $\pi$. We can apply the same reasoning to deduce that for $\mathcal{B}\rightarrow 0$, the unstable equilibria are such that $\theta^{\ast}(\pm 1) \rightarrow (k+\frac{1}{2})\pi$, as previously claimed before Section \ref{sec:G0}. For weaker anchoring, the director profile has greater freedom to reorient at the boundaries and escape from the energetically expensive fixed rotation imposed by large winding numbers. 
For $\mathcal{G} = 0$, $\mathcal{B}_{i}^{\ast}= \mathcal{B}_{-i}^{\ast}$ ($i\in \mathds{N}$). For $\mathcal{B} > \mathcal{B}_{1}^{\ast}$, $\theta_{a_{0}}^{\ast} $ and $\theta_{\tilde{a}_{0}}^{\ast}$ are the only constant steady states of system \eqref{systemstatic}. For simplicity,
in what follows we denote the equilibrium solutions as $\theta^{\ast}_{a}$, where $\theta^{\ast}_{a}= \theta^{\ast}_{a_{n}}$ if it is of Type I and $\theta^{\ast}_{a}= \theta^{\ast}_{\tilde{a}_{n}}$ if it is of Type  II.     
  
\subsection{Fluid flow ($\mathcal{G}>0$)}\label{sec:gpositive}
Next, we  study the static equilibria of the system (\ref{systemstatic})
when we apply a pressure difference $\mathcal{G}>0$ across the
microfluidic channel, inducing a fluid flow. The solutions are computed
numerically for all values of $\mathcal{G}$ using \textit{Chebfun} via the method of continuation.\cite{allgower2003} When the $\mathcal{G}=0$ solution $\theta^{\ast}_{a}$ is taken as the initial condition (see \textsection\ref{sec:G0}), the corresponding solution with $\mathcal{G}>0$ is denoted by $\theta^{\ast}_{a,\mathcal{G}}$.

\subsubsection{Asymptotics when $\mathcal{G} \ll 1$}\label{sec:Gsmall}
When $\mathcal{G}\ll 1$, we can approximate
$\theta_{a,\mathcal{G}}^{\ast}$ by the expansion \linebreak
$\theta_{a,\mathcal{G}}^{\ast}(z)= \theta_{a}^{\ast}(z) +
\mathcal{G}\theta_{\mathcal{G}}^{(1)}(z) +\cdots$, where $\theta_{a}^\ast$
is the corresponding solution for $\mathcal{G}=0$. It is
straightforward to verify that $\theta_{\mathcal{G}}^{(1)}$
satisfies
\begin{align}
\label{eq:theta 1 G>0}
\left\{
\begin{array}{l} 
\displaystyle \frac{{\rm d}^{2} \theta^{(1)}_{\mathcal{G}}}{{\rm d} z^{2}}(z)  = z Q({\theta_{a}^\ast}(z)) \quad z\in (-1,1)\\
\mathcal{B} \displaystyle \frac{{\rm d} \theta^{(1)}_{\mathcal{G}}}{{\rm d}z}(1)  =
-2{\theta^{(1)}_{\mathcal{G}}}(1) \cos(2{\theta^\ast_{a}}(1)),   \\ 
\mathcal{B} \displaystyle \frac{{\rm d} \theta^{(1)}_{\mathcal{G}}}{{\rm d}z}(-1)  =
2{\theta^{(1)}_{\mathcal{G}}}(-1) \cos(2{\theta^\ast_{a}}(-1)),  \end{array}\right.
\end{align}
where $Q(s)=-{m(s)}/{g(s)}$.
The solution to \eqref{eq:theta 1 G>0} is given by
\begin{equation}\label{eq:new}
\theta_{\mathcal{G}}^{(1)}(z)= J(z) + Cz + D,
\end{equation}
where
\begin{align}
&I(r)=
\int_{0}^{r} s Q(as+b){\rm d}s, \hspace{20mm} J(z)=
\int_{0}^{z} I(r){\rm d} r, \\
&C   = \displaystyle \frac{2(-1)^{k}\cos(2a)\big( J(-1)-J(1)\big)-\mathcal{B}\big(I(1)+I(-1)\big)}{2\mathcal{B}+4(-1)^{k}\cos(2a)}, \\
&D= \displaystyle -\frac{1}{2}\big( J(1)+J(-1) \big)+ \frac{\mathcal{B}(-1)^{k}\big(I(-1)-I(1)\big)}{4\cos(2a)},
\end{align}
with $b=k=0$ for Type I solutions where $a$ satisfies \eqref{eq:type1m0} and $b=\frac{\pi}{2}$ and $k=1$ for Type II solutions, where $a$ satisfies \eqref{eq:type2m0}.
\begin{figure}[ht!]
\centering
\subfigure[Comparison of asymptotic solution given by \eqref{eq:new} (dashed) with the full numerical solution to \eqref{systemstatic} (solid)]{
%
%
\begin{psfrags}%
\psfragscanon%
%
\psfrag{s03}[t][t]{\color[rgb]{0,0,0}\setlength{\tabcolsep}{0pt}\begin{tabular}{c}{\LARGE $z$}\end{tabular}}%
\psfrag{s04}[b][b]{\color[rgb]{0,0,0}\setlength{\tabcolsep}{0pt}\begin{tabular}{c}{\LARGE $\theta^{\ast}_{a_{0},\mathcal{G}}(z) $}\end{tabular}}%
\psfrag{s05}[l][l]{\color[rgb]{0,0,0}\setlength{\tabcolsep}{0pt}\begin{tabular}{l}{\Large $\mathcal{G}=$0.1}\end{tabular}}%
\psfrag{s06}[l][l]{\color[rgb]{0,0,0}\setlength{\tabcolsep}{0pt}\begin{tabular}{l}{\Large $\mathcal{G}=$2}\end{tabular}}%
\psfrag{s07}[l][l]{\color[rgb]{0,0,0}\setlength{\tabcolsep}{0pt}\begin{tabular}{l}{\Large $\mathcal{G}=$5}\end{tabular}}%
\psfrag{s08}[l][l]{\color[rgb]{0,0,0}\setlength{\tabcolsep}{0pt}\begin{tabular}{l}{\Large $\mathcal{G}=$10}\end{tabular}}%
%
\psfrag{x01}[t][t]{-1}%
\psfrag{x02}[t][t]{-0.8}%
\psfrag{x03}[t][t]{-0.6}%
\psfrag{x04}[t][t]{-0.4}%
\psfrag{x05}[t][t]{-0.2}%
\psfrag{x06}[t][t]{0}%
\psfrag{x07}[t][t]{0.2}%
\psfrag{x08}[t][t]{0.4}%
\psfrag{x09}[t][t]{0.6}%
\psfrag{x10}[t][t]{0.8}%
\psfrag{x11}[t][t]{1}%
%
\psfrag{v01}[r][r]{-0.8}%
\psfrag{v02}[r][r]{-0.6}%
\psfrag{v03}[r][r]{-0.4}%
\psfrag{v04}[r][r]{-0.2}%
\psfrag{v05}[r][r]{0}%
\psfrag{v06}[r][r]{0.2}%
\psfrag{v07}[r][r]{0.4}%
\psfrag{v08}[r][r]{0.6}%
\psfrag{v09}[r][r]{0.8}%
%
\resizebox{8cm}{!}{\includegraphics{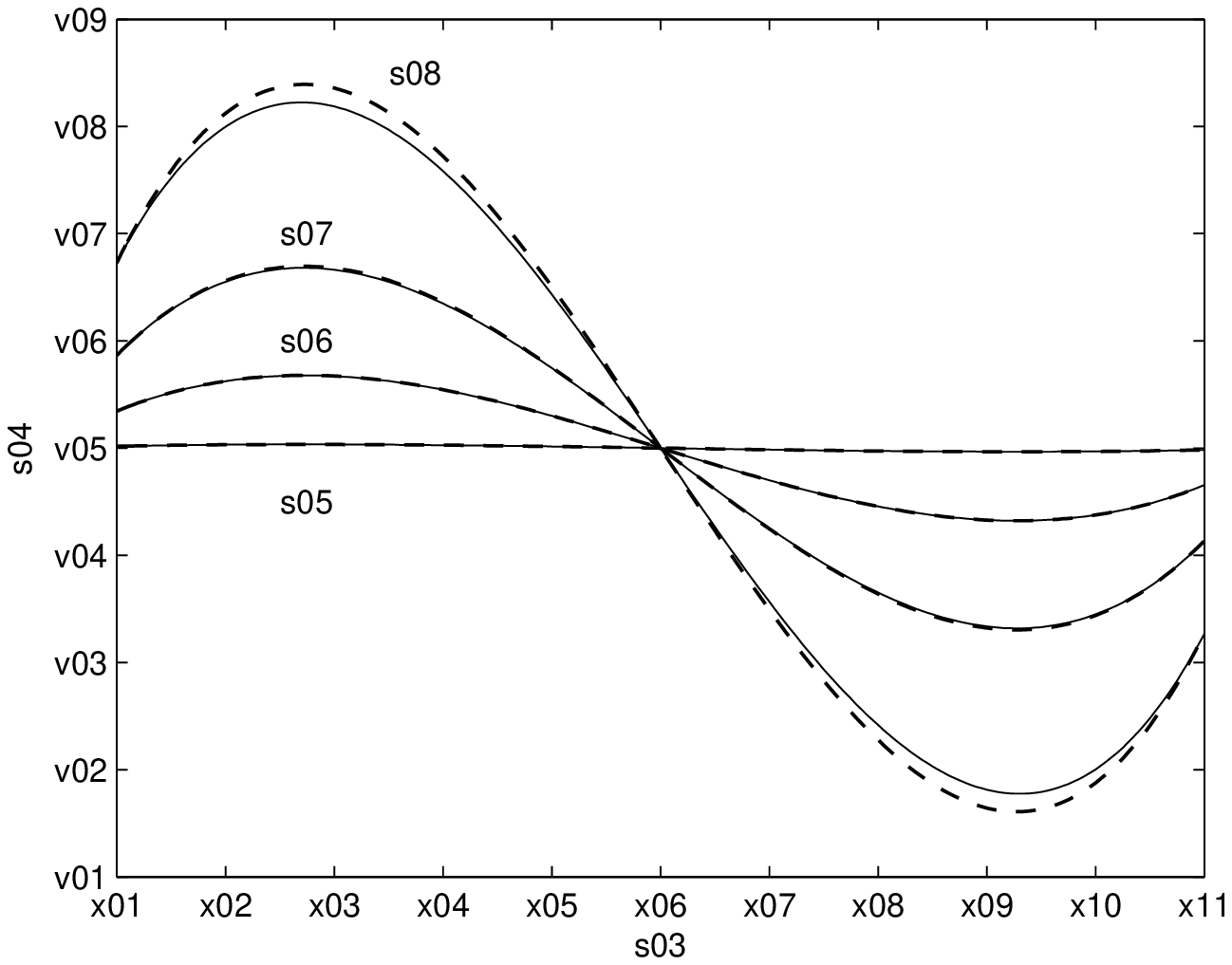}}%
\end{psfrags}%
%
}
\subfigure[$\ve{n}$ with $\theta^{\ast}_{a_{0},0.5}$]{
%
%
\begin{psfrags}%
\psfragscanon%
%
\psfrag{s03}[t][t]{\color[rgb]{0,0,0}\setlength{\tabcolsep}{0pt}\begin{tabular}{c}{\Large $x$}\end{tabular}}%
\psfrag{s04}[b][b]{\color[rgb]{0,0,0}\setlength{\tabcolsep}{0pt}\begin{tabular}{c}{\Large $z$}\end{tabular}}%
%
\psfrag{x01}[t][t]{-0.005}%
\psfrag{x02}[t][t]{0.005}%
%
\psfrag{v01}[r][r]{-1}%
\psfrag{v02}[r][r]{1}%
%
\resizebox{2.8cm}{!}{\includegraphics{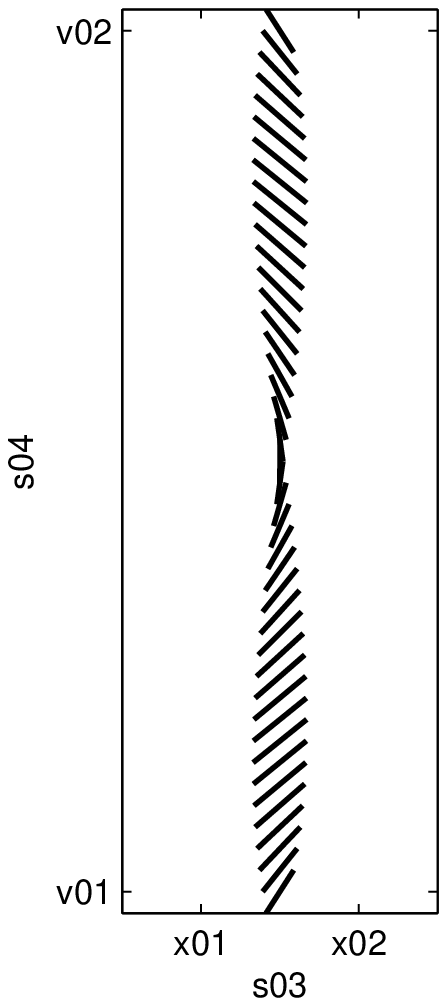}}%
\end{psfrags}%
%
}
\caption{Static equilibra $\theta^{\ast}_{a_{0},\mathcal{G}}$ when $\mathcal{B}=\frac{1}{3}$. (b) We note that $\ve{n} \approx (0,0,1)$ but different scales have been used in the $x$ and $z$ axis to allow the reader to appreciate the change between $\theta^{\ast}_{a_{0},0.5}$ and $\theta^{\ast}_{a_{0},0}$ (corresponding to $\ve{n}=(0,0,1)$).}
\label{fig:asymptoticsnearg0t0}
\end{figure}
\begin{figure}[ht!]
\centering
\subfigure[Comparison of asymptotic solution given by \eqref{eq:new} (dashed) with the full numerical solution to \eqref{systemstatic} (solid)]{
%
%
\begin{psfrags}%
\psfragscanon%
%
\psfrag{s03}[t][t]{\color[rgb]{0,0,0}\setlength{\tabcolsep}{0pt}\begin{tabular}{c}{\LARGE $z$}\end{tabular}}%
\psfrag{s04}[b][b]{\color[rgb]{0,0,0}\setlength{\tabcolsep}{0pt}\begin{tabular}{c}{\LARGE $\theta^{\ast}_{\tilde{a}_{1},\mathcal{G}}(z) $}\end{tabular}}%
\psfrag{s05}[l][l]{\color[rgb]{0,0,0}\setlength{\tabcolsep}{0pt}\begin{tabular}{l}{\Large $\mathcal{G}=$7}\end{tabular}}%
\psfrag{s06}[l][l]{\color[rgb]{0,0,0}\setlength{\tabcolsep}{0pt}\begin{tabular}{l}{\Large  $\mathcal{G}=$0.1}\end{tabular}}%
\psfrag{s07}[l][l]{\color[rgb]{0,0,0}\setlength{\tabcolsep}{0pt}\begin{tabular}{l}{\Large $\mathcal{G}=$2}\end{tabular}}%
%
\psfrag{x01}[t][t]{0}%
\psfrag{x02}[t][t]{0.1}%
\psfrag{x03}[t][t]{0.2}%
\psfrag{x04}[t][t]{0.3}%
\psfrag{x05}[t][t]{0.4}%
\psfrag{x06}[t][t]{0.5}%
\psfrag{x07}[t][t]{0.6}%
\psfrag{x08}[t][t]{0.7}%
\psfrag{x09}[t][t]{0.8}%
\psfrag{x10}[t][t]{0.9}%
\psfrag{x11}[t][t]{1}%
\psfrag{x12}[t][t]{-1}%
\psfrag{x13}[t][t]{-0.8}%
\psfrag{x14}[t][t]{-0.6}%
\psfrag{x15}[t][t]{-0.4}%
\psfrag{x16}[t][t]{-0.2}%
\psfrag{x17}[t][t]{0}%
\psfrag{x18}[t][t]{0.2}%
\psfrag{x19}[t][t]{0.4}%
\psfrag{x20}[t][t]{0.6}%
\psfrag{x21}[t][t]{0.8}%
\psfrag{x22}[t][t]{1}%
%
\psfrag{v01}[r][r]{0}%
\psfrag{v02}[r][r]{0.1}%
\psfrag{v03}[r][r]{0.2}%
\psfrag{v04}[r][r]{0.3}%
\psfrag{v05}[r][r]{0.4}%
\psfrag{v06}[r][r]{0.5}%
\psfrag{v07}[r][r]{0.6}%
\psfrag{v08}[r][r]{0.7}%
\psfrag{v09}[r][r]{0.8}%
\psfrag{v10}[r][r]{0.9}%
\psfrag{v11}[r][r]{1}%
\psfrag{v12}[r][r]{0}%
\psfrag{v13}[r][r]{0.5}%
\psfrag{v14}[r][r]{1}%
\psfrag{v15}[r][r]{1.5}%
\psfrag{v16}[r][r]{2}%
\psfrag{v17}[r][r]{2.5}%
\psfrag{v18}[r][r]{3}%
%
\resizebox{8cm}{!}{\includegraphics{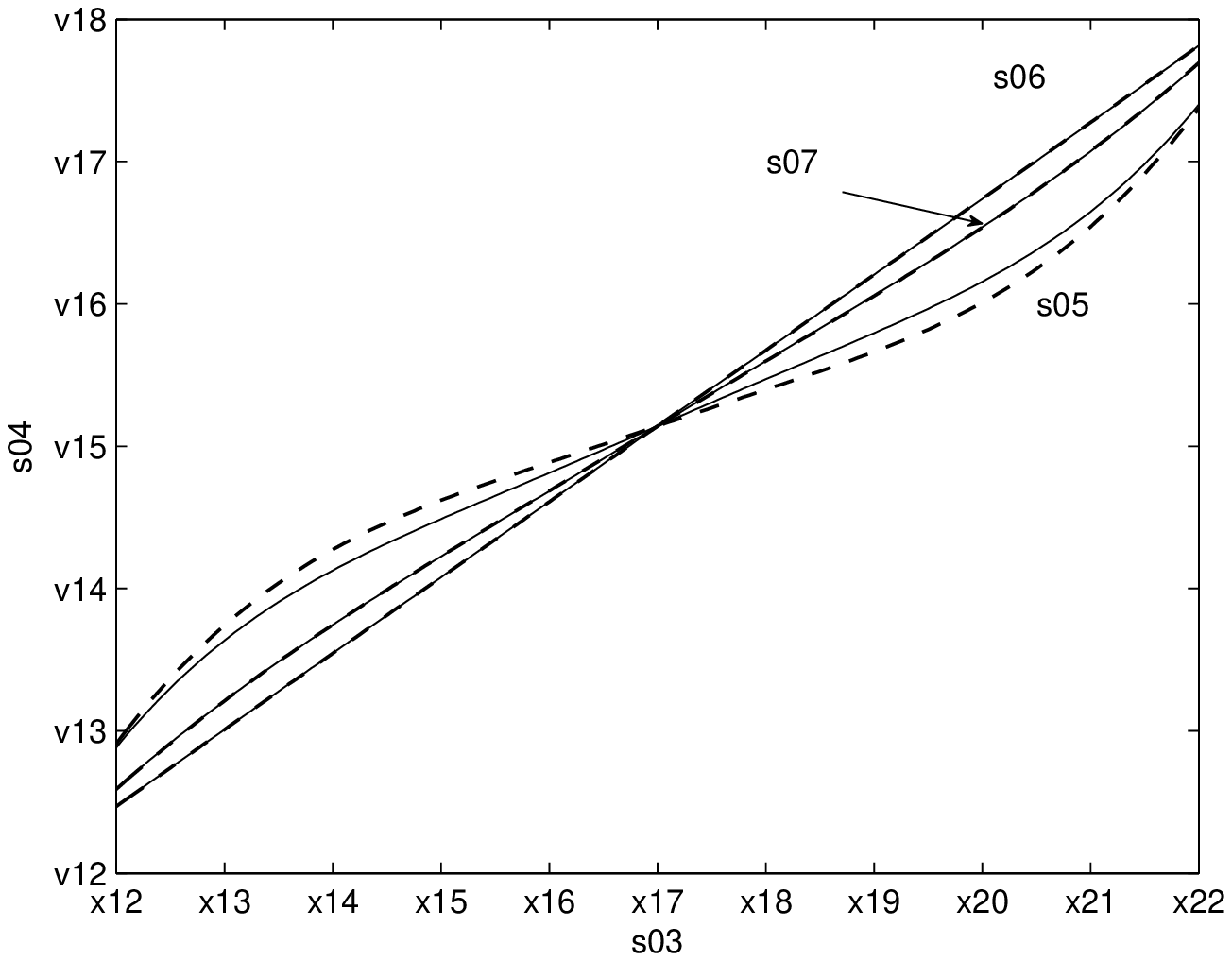}}%
\end{psfrags}%
%
}
\subfigure[$\ve{n}$ with $\theta^{\ast}_{\tilde{a}_{1},0.5}$]{
%
%
\begin{psfrags}%
\psfragscanon%
%
\psfrag{s03}[t][t]{\color[rgb]{0,0,0}\setlength{\tabcolsep}{0pt}\begin{tabular}{c}{\Large $x$}\end{tabular}}%
\psfrag{s04}[b][b]{\color[rgb]{0,0,0}\setlength{\tabcolsep}{0pt}\begin{tabular}{c}{\Large $z$}\end{tabular}}%
%
\psfrag{x01}[t][t]{-0.1}%
\psfrag{x02}[t][t]{0}%
\psfrag{x03}[t][t]{0.1}%
%
\psfrag{v01}[r][r]{-1}%
\psfrag{v02}[r][r]{1}%
%
\resizebox{2.8cm}{!}{\includegraphics{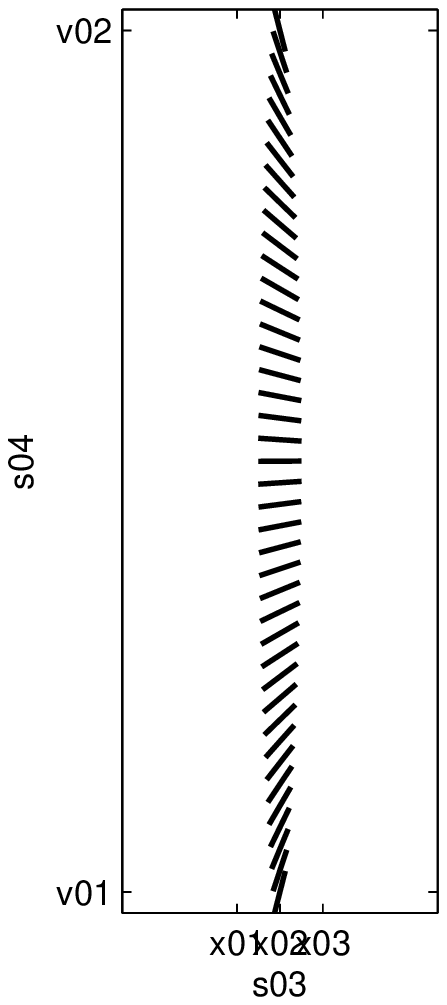}}%
\end{psfrags}%
%
}
\caption{Static equilibria $\theta^{\ast}_{\tilde{a}_{1},\mathcal{G}}$ when $\mathcal{B}=\frac{1}{3}$. (b)  In contrast with Figure \ref{fig:asymptoticsnearg0t0}(b), here $x$ and $z$ axis have the same scale, which corresponds to the real configuration of the molecules.}
\label{fig:asymptoticsnearg0t2}
\end{figure}
We validate the asymptotic analysis performed above by numerically computing the equilibria $\theta^{\ast}_{a, \mathcal{G}}$ of (\ref{systemtdep}) for small values of
$\mathcal{G}$ by solving (\ref{systemstatic}) with \textit{Chebfun} and comparing this with the asymptotic result
\eqref{eq:new}. When $\theta_{a}^{\ast}=\theta^{\ast}_{a_{0}} \equiv 0$ and $\theta_{a}^{\ast}=\theta^{\ast}_{\tilde{a}_{1}}$ the asymptotic solution approximates the actual solution for values of $\mathcal{G}$ significantly
beyond the expected regime (see respectively Figures \ref{fig:asymptoticsnearg0t0}(a) and \ref{fig:asymptoticsnearg0t2}(a), where we find that the asymptotic
solution approximates the full numerical solution well for values of $\mathcal{G}$ as large as 7). Figures \ref{fig:asymptoticsnearg0t0}(b) and \ref{fig:asymptoticsnearg0t2}(b) show the director field $\ve{n}$ associated with the equilibria $\theta_{a_{0},\mathcal{G}}^{\ast}$ and $\theta_{\tilde{a}_{1},\mathcal{G}}^{\ast}$, computed when $\mathcal{G}=0.5$ and $\mathcal{B}=\frac{1}{3}$. We chose a moderate anchoring strength to illustrate the differences between the numerics and asymptotics clearly. The asymptotic approximations rapidly improve as $\mathcal{B}\to 0$.
\subsubsection{Asymptotics when $\mathcal{G} \rightarrow \infty$}\label{sec:Gbig}
For ${\mathcal{G}}\gg 1$, we can perform a similar asymptotic expansion of the form
$\theta^{\ast}_{\mathcal{G}}(z)= \theta^{(0)}_{\mathcal{G}}(z) + (1/\mathcal{G})\theta^{(1)}_{\mathcal{G}}(z) + \cdots$. Substituting this expansion into (\ref{systemstatic}) and equating terms at leading order gives
\begin{subequations}
\label{outerGinfty}
\begin{align}
z Q(\theta^{(0)}_{\mathcal{G}}(z))&=0, & \,\,\,z\in (-1,1)  \label{outerGinfty1}
 \\
\mathcal{B}\displaystyle \frac{{\rm d}\theta^{(0)}_{\mathcal{G}}}{{\rm d}z}(1)&=-\sin(2\theta_{\mathcal{G}}^{(0)}(1)),
 &   \label{outerGinfty2} \\
\mathcal{B}\displaystyle \frac{{\rm d}\theta^{(0)}_{\mathcal{G}}}{{\rm d}z}(-1)&=\sin(2\theta_{\mathcal{G}}^{(0)}(-1)). & \label{outerGinfty3}
\end{align}
\end{subequations}
Equation (\ref{outerGinfty1}) implies that $\theta^{(0)}_{\mathcal{G}}(0)$ can take arbitrary values in $\mathds{R}$ and
\begin{align}
\label{solouterGinfty}
\theta_{\mathcal{G}}^{(0)}(z) \equiv \pm \arctan\lr{ \sqrt{\frac{\alpha_{2}}{\alpha_{3}}}} + k\pi\equiv \sigma^\pm_k \hspace*{1cm}\forall z \neq 0,
\end{align}
with $k\in \mathds{Z}$ arbitrary. However, the boundary conditions (\ref{outerGinfty}\textit{b.c}) are not satisfied by \eqref{solouterGinfty} and hence we expect to find boundary layers near $z=-1$, 0 and 1, in order to match the boundary conditions. The solution in the two outer regions $-1 < z < 0$ and $0 < z < 1$ are given by \eqref{solouterGinfty} for any two particular integer values of $k$, say $k_1$ and $k_2$.

Near $z=-1$, we rescale in (\ref{systemstatic}) by introducing the variable $\eta= {\sqrt{\mathcal{G}}}({z+1})$ and perform an asymptotic expansion in powers of $1 /\sqrt{\mathcal{G}}$. The corresponding leading-order term in $\mathcal{G}$, $\theta_{\textrm{L}, \mathcal{G}}^{(0)}(\eta)$, is a solution of
\begin{subequations}
\label{inner1Ginfty}
\begin{align}
\frac{{\rm d}^{2} \theta^{(0)}_{\textrm{L},\mathcal{G}}}{{\rm d} \eta^{2}}(\eta)= - Q(\theta^{(0)}_{\textrm{L},\mathcal{G}}(\eta)), & \eta>0  \label{inner1Ginfty1}\\
\bar{\mathcal{B}} \frac{{\rm d} \theta^{(0)}_{\rm{L},\mathcal{G}}}{{\rm d} \eta}(0)= \sin(2\theta^{(0)}_{\textrm{L},\mathcal{G}}(0)), & \label{inner1Ginfty2}\\
\lim_{\eta \rightarrow \infty} \theta^{(0)}_{\textrm{L},\mathcal{G}}(\eta)= \sigma_{k_1}^\pm, & \label{inner1Ginfty3}
\end{align}
\end{subequations}
where we have rescaled $\bar{\mathcal{B}}= \sqrt{\mathcal{G}}\mathcal{B}$ assuming that $\bar{\mathcal{B}}=O(1)$ to obtain the richest asymptotic limit. We point out that the asymptotic analysis could be done without this assumption. Then \eqref{inner1Ginfty2} would be $\mathcal{B} \frac{{\rm d}\theta^{(0)}_{\rm{L},\mathcal{G}}}{{\rm d} \eta} =0$ and $ \theta^{(0)}_{\textrm{L},\mathcal{G}}(\eta)= \sigma_{k_{1}}^{\pm}$. We would need to use the second term, $\theta^{(1)}_{\textrm{L},\mathcal{G}}$, of the asymptotic expansion (at least) and the results with these two terms would be worse than those obtained here. Equation (\ref{inner1Ginfty}\textit{c}) is the matching condition between $\theta^{(0)}_{\textrm{L},\mathcal{G}}$ and $\theta^{(0)}_{\mathcal{G}}$.

Near $z=0$, we set $\xi= \mathcal{G}^{1/3}z$ and the corresponding leading-order term, $\theta^{(0)}_{\textrm{C},\mathcal{G}}(\xi)$, satisfies 
\begin{subequations}
\label{inner2Ginfty}
\begin{align}
\frac{{\rm d}^{2} \theta^{(0)}_{\textrm{C},  \mathcal{G} }}{{\rm d}\xi^{2}}= \xi Q(\theta^{(0)}_{\textrm{C},\mathcal{G}}(\xi)),
 & \,\,\,\xi\in (-\infty,\infty),  \\
\lim_{\xi \rightarrow -\infty} \theta^{(0)}_{\textrm{C},\mathcal{G}}(\xi)= \sigma_{k_1}^{\pm} , &\\
\lim_{\xi \rightarrow \infty} \theta^{(0)}_{\textrm{C},\mathcal{G}}(\eta)= \sigma_{k_2}^{\pm}, &
\end{align}
\end{subequations}
where (\ref{inner2Ginfty}\textit{b,c}) describe the matching conditions.

Finally, we introduce the variable $\zeta= \sqrt{\mathcal{G}}(1-z)$  near $z=1$ and $\theta^{(0)}_{\textrm{R},\mathcal{G}}(\zeta)$, the leading--order solution in $\mathcal{G}$, satisfies
\begin{subequations}
\label{inner3Ginfty}
\begin{align}
\frac{{\rm d}^{2}\theta^{(0)}_{\textrm{R},\mathcal{G}}}{{\rm d} \zeta^{2}}= Q(\theta^{(0)}_{\textrm{R},\mathcal{G}}(\zeta)), & \,\,\,\zeta>0,  \\
\bar{\mathcal{B}} \frac{{\rm d} \theta^{(0)}_{\rm{L},\mathcal{G}}}{{\rm d} \eta}(0)= \sin(2\theta^{(0)}_{\textrm{L},\mathcal{G}}(0)),&
\\
\lim_{\zeta\rightarrow \infty} \theta^{(0)}_{\textrm{R},\mathcal{G}}(\zeta) = \sigma_{k_2}^{\pm},
\end{align}
\end{subequations}
where (\ref{inner3Ginfty}\textit{c}) is the matching condition.

\begin{figure}[ht!]
\centering
\subfigure[Comparison of asymptotic solution given by \eqref{eq:new} (dashed) with the full numerical solution to \eqref{systemstatic} (solid)]{
%
%
\begin{psfrags}%
\psfragscanon%
%
\psfrag{s03}[t][t]{\color[rgb]{0,0,0}\setlength{\tabcolsep}{0pt}\begin{tabular}{c}{\LARGE $z$}\end{tabular}}%
\psfrag{s04}[b][b]{\color[rgb]{0,0,0}\setlength{\tabcolsep}{0pt}\begin{tabular}{c}{\LARGE $\theta^{\ast}_{a_{0},\mathcal{G}}(z) $}\end{tabular}}%
\psfrag{s05}[l][l]{\color[rgb]{0,0,0}\setlength{\tabcolsep}{0pt}\begin{tabular}{l}{\Large increasing $\mathcal{G}$}\end{tabular}}%
\psfrag{s06}[l][l]{\color[rgb]{0,0,0}\setlength{\tabcolsep}{0pt}\begin{tabular}{l}$\mathcal{G}=$1000 (asym)\end{tabular}}%
\psfrag{s07}[l][l]{\color[rgb]{0,0,0}\setlength{\tabcolsep}{0pt}\begin{tabular}{l}$\mathcal{G}=$1000\end{tabular}}%
\psfrag{s08}[l][l]{\color[rgb]{0,0,0}\setlength{\tabcolsep}{0pt}\begin{tabular}{l}$\mathcal{G}=$500 (asym)\end{tabular}}%
\psfrag{s09}[l][l]{\color[rgb]{0,0,0}\setlength{\tabcolsep}{0pt}\begin{tabular}{l}$\mathcal{G}=$500\end{tabular}}%
\psfrag{s10}[l][l]{\color[rgb]{0,0,0}\setlength{\tabcolsep}{0pt}\begin{tabular}{l}$\mathcal{G}=$100 (asym)\end{tabular}}%
\psfrag{s11}[l][l]{\color[rgb]{0,0,0}\setlength{\tabcolsep}{0pt}\begin{tabular}{l}$\mathcal{G}=$100\end{tabular}}%
%
\psfrag{x01}[t][t]{0}%
\psfrag{x02}[t][t]{0.1}%
\psfrag{x03}[t][t]{0.2}%
\psfrag{x04}[t][t]{0.3}%
\psfrag{x05}[t][t]{0.4}%
\psfrag{x06}[t][t]{0.5}%
\psfrag{x07}[t][t]{0.6}%
\psfrag{x08}[t][t]{0.7}%
\psfrag{x09}[t][t]{0.8}%
\psfrag{x10}[t][t]{0.9}%
\psfrag{x11}[t][t]{1}%
\psfrag{x12}[t][t]{-1}%
\psfrag{x13}[t][t]{-0.8}%
\psfrag{x14}[t][t]{-0.6}%
\psfrag{x15}[t][t]{-0.4}%
\psfrag{x16}[t][t]{-0.2}%
\psfrag{x17}[t][t]{0}%
\psfrag{x18}[t][t]{0.2}%
\psfrag{x19}[t][t]{0.4}%
\psfrag{x20}[t][t]{0.6}%
\psfrag{x21}[t][t]{0.8}%
\psfrag{x22}[t][t]{1}%
%
\psfrag{v01}[r][r]{0}%
\psfrag{v02}[r][r]{0.1}%
\psfrag{v03}[r][r]{0.2}%
\psfrag{v04}[r][r]{0.3}%
\psfrag{v05}[r][r]{0.4}%
\psfrag{v06}[r][r]{0.5}%
\psfrag{v07}[r][r]{0.6}%
\psfrag{v08}[r][r]{0.7}%
\psfrag{v09}[r][r]{0.8}%
\psfrag{v10}[r][r]{0.9}%
\psfrag{v11}[r][r]{1}%
\psfrag{v12}[r][r]{-1.5}%
\psfrag{v13}[r][r]{-1}%
\psfrag{v14}[r][r]{-0.5}%
\psfrag{v15}[r][r]{0}%
\psfrag{v16}[r][r]{0.5}%
\psfrag{v17}[r][r]{1}%
\psfrag{v18}[r][r]{1.5}%
%
\resizebox{8cm}{!}{\includegraphics{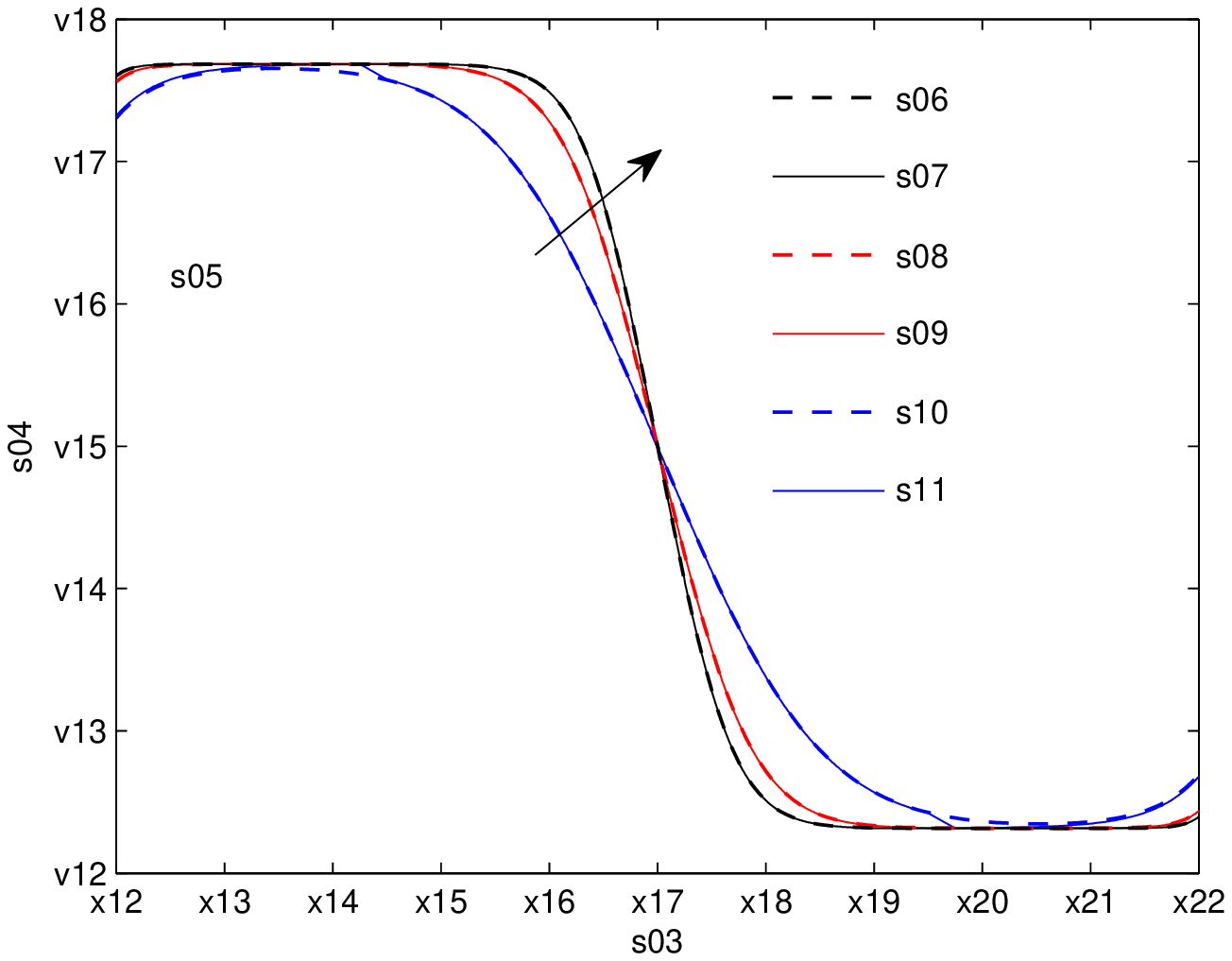}}%
\end{psfrags}%
%
}
\subfigure[$\ve{n}$ with $\theta^{\ast}_{a_{0},100}$]{
%
%
\begin{psfrags}%
\psfragscanon%
%
\psfrag{s03}[t][t]{\color[rgb]{0,0,0}\setlength{\tabcolsep}{0pt}\begin{tabular}{c}{\Large $x$}\end{tabular}}%
\psfrag{s04}[b][b]{\color[rgb]{0,0,0}\setlength{\tabcolsep}{0pt}\begin{tabular}{c}{\Large $z$}\end{tabular}}%
%
\psfrag{x01}[t][t]{-0.1}%
\psfrag{x02}[t][t]{0}%
\psfrag{x03}[t][t]{0.1}%
%
\psfrag{v01}[r][r]{-1}%
\psfrag{v02}[r][r]{0}%
\psfrag{v03}[r][r]{1}%
%
\resizebox{2.8cm}{!}{\includegraphics{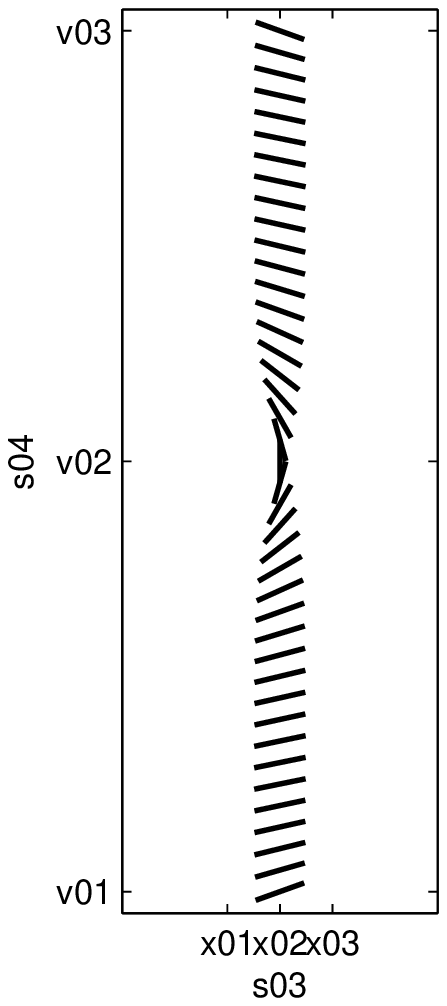}}%
\end{psfrags}%
%
}
\caption{Static equilibria $\theta^{\ast}_{a_{0},\mathcal{G}}$ with $\mathcal{G}\rightarrow \infty$ and $\mathcal{B}=\frac{1}{3}$.}
\label{fig:asymptoticsnearginft0}
\end{figure}

\begin{figure}[ht!]
\centering
\subfigure[Comparison of asymptotic solution given by \eqref{eq:new} (dashed) with the full numerical solution to \eqref{systemstatic} (solid)]{
%
%
\begin{psfrags}%
\psfragscanon%
%
\psfrag{s03}[t][t]{\color[rgb]{0,0,0}\setlength{\tabcolsep}{0pt}\begin{tabular}{c}{\LARGE $z$}\end{tabular}}%
\psfrag{s04}[b][b]{\color[rgb]{0,0,0}\setlength{\tabcolsep}{0pt}\begin{tabular}{c}{\LARGE $\theta_{\tilde{a}_{1},\mathcal{G}}^{\ast}(z) $}\end{tabular}}%
\psfrag{s05}[l][l]{\color[rgb]{0,0,0}\setlength{\tabcolsep}{0pt}\begin{tabular}{l}{\Large increasing $\mathcal{G}$}\end{tabular}}%
\psfrag{s06}[l][l]{\color[rgb]{0,0,0}\setlength{\tabcolsep}{0pt}\begin{tabular}{l}$\mathcal{G}=$1000 (asym)\end{tabular}}%
\psfrag{s07}[l][l]{\color[rgb]{0,0,0}\setlength{\tabcolsep}{0pt}\begin{tabular}{l}$\mathcal{G}=$1000\end{tabular}}%
\psfrag{s08}[l][l]{\color[rgb]{0,0,0}\setlength{\tabcolsep}{0pt}\begin{tabular}{l}$\mathcal{G}=$500 (asym)\end{tabular}}%
\psfrag{s09}[l][l]{\color[rgb]{0,0,0}\setlength{\tabcolsep}{0pt}\begin{tabular}{l}$\mathcal{G}=$500\end{tabular}}%
\psfrag{s10}[l][l]{\color[rgb]{0,0,0}\setlength{\tabcolsep}{0pt}\begin{tabular}{l}$\mathcal{G}=$100 (asym)\end{tabular}}%
\psfrag{s11}[l][l]{\color[rgb]{0,0,0}\setlength{\tabcolsep}{0pt}\begin{tabular}{l}$\mathcal{G}=$100\end{tabular}}%
%
\psfrag{x01}[t][t]{0}%
\psfrag{x02}[t][t]{0.1}%
\psfrag{x03}[t][t]{0.2}%
\psfrag{x04}[t][t]{0.3}%
\psfrag{x05}[t][t]{0.4}%
\psfrag{x06}[t][t]{0.5}%
\psfrag{x07}[t][t]{0.6}%
\psfrag{x08}[t][t]{0.7}%
\psfrag{x09}[t][t]{0.8}%
\psfrag{x10}[t][t]{0.9}%
\psfrag{x11}[t][t]{1}%
\psfrag{x12}[t][t]{-1}%
\psfrag{x13}[t][t]{-0.8}%
\psfrag{x14}[t][t]{-0.6}%
\psfrag{x15}[t][t]{-0.4}%
\psfrag{x16}[t][t]{-0.2}%
\psfrag{x17}[t][t]{0}%
\psfrag{x18}[t][t]{0.2}%
\psfrag{x19}[t][t]{0.4}%
\psfrag{x20}[t][t]{0.6}%
\psfrag{x21}[t][t]{0.8}%
\psfrag{x22}[t][t]{1}%
%
\psfrag{v01}[r][r]{0}%
\psfrag{v02}[r][r]{0.1}%
\psfrag{v03}[r][r]{0.2}%
\psfrag{v04}[r][r]{0.3}%
\psfrag{v05}[r][r]{0.4}%
\psfrag{v06}[r][r]{0.5}%
\psfrag{v07}[r][r]{0.6}%
\psfrag{v08}[r][r]{0.7}%
\psfrag{v09}[r][r]{0.8}%
\psfrag{v10}[r][r]{0.9}%
\psfrag{v11}[r][r]{1}%
\psfrag{v12}[r][r]{1.1}%
\psfrag{v13}[r][r]{1.2}%
\psfrag{v14}[r][r]{1.3}%
\psfrag{v15}[r][r]{1.4}%
\psfrag{v16}[r][r]{1.5}%
\psfrag{v17}[r][r]{1.6}%
\psfrag{v18}[r][r]{1.7}%
\psfrag{v19}[r][r]{1.8}%
\psfrag{v20}[r][r]{1.9}%
\psfrag{v21}[r][r]{2}%
\psfrag{v22}[r][r]{2.1}%
%
\resizebox{8cm}{!}{\includegraphics{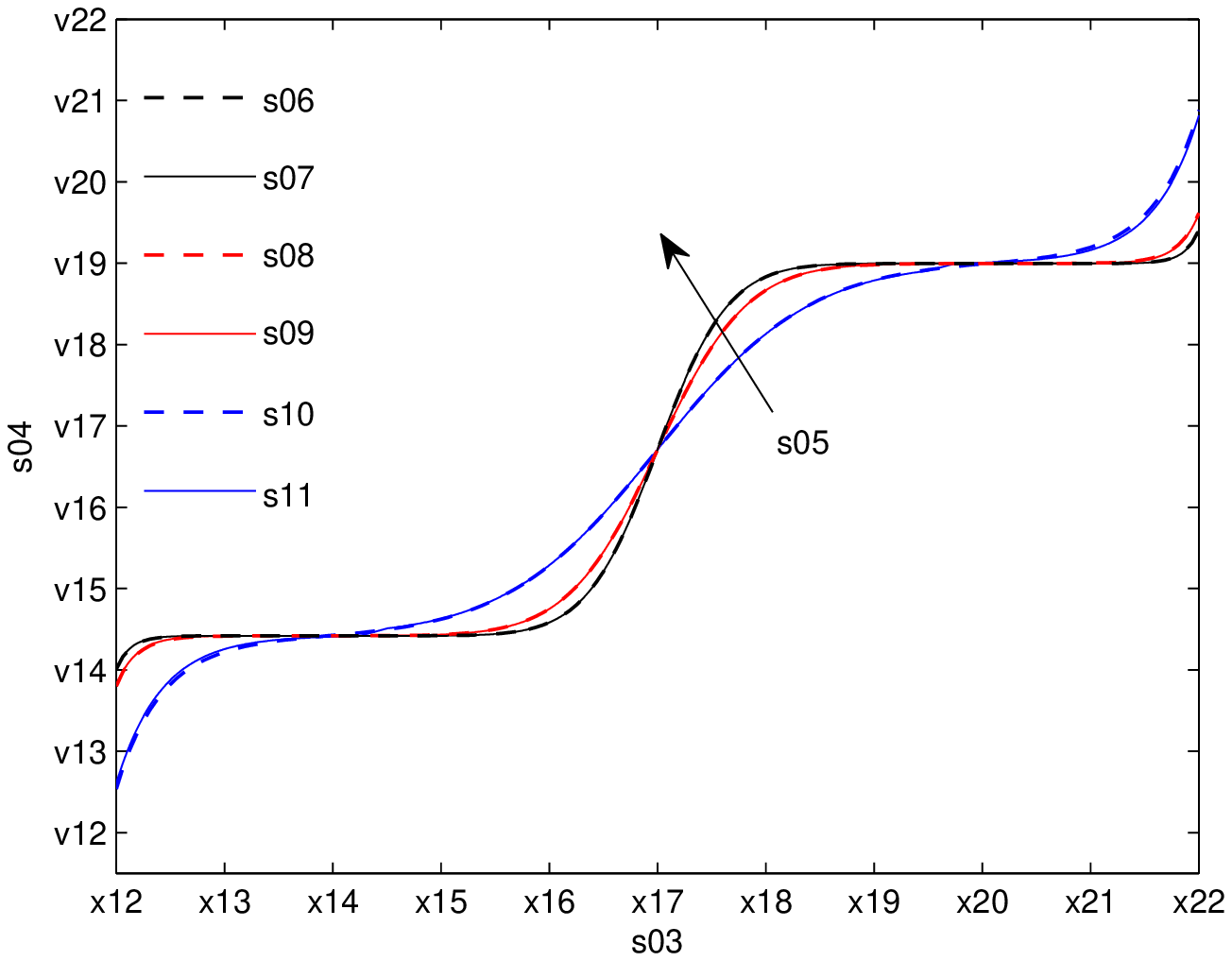}}%
\end{psfrags}%
%
}
\subfigure[$\ve{n}$ with $\theta^{\ast}_{\tilde{a}_{1},100}$]{
%
%
\begin{psfrags}%
\psfragscanon%
%
\psfrag{s03}[t][t]{\color[rgb]{0,0,0}\setlength{\tabcolsep}{0pt}\begin{tabular}{c}{\Large $x$}\end{tabular}}%
\psfrag{s04}[b][b]{\color[rgb]{0,0,0}\setlength{\tabcolsep}{0pt}\begin{tabular}{c}{\Large $z$}\end{tabular}}%
%
\psfrag{x01}[t][t]{-0.1}%
\psfrag{x02}[t][t]{0}%
\psfrag{x03}[t][t]{0.1}%
%
\psfrag{v01}[r][r]{-1}%
\psfrag{v02}[r][r]{0}%
\psfrag{v03}[r][r]{1}%
%
\resizebox{2.8cm}{!}{\includegraphics{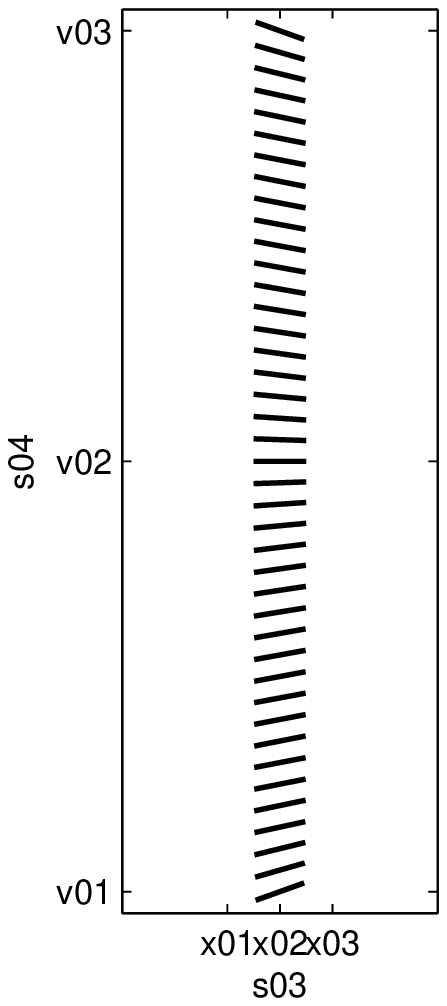}}%
\end{psfrags}%
%
}
\caption{Static equilibria $\theta^{\ast}_{\tilde{a}_{1},\mathcal{G}}$ when $\mathcal{G}\rightarrow \infty$ and $\mathcal{B}=\frac{1}{3}$.}
\label{fig:asymptoticsnearginft2}
\end{figure}

We numerically solve the three boundary layer problems (\ref{inner1Ginfty}), (\ref{inner2Ginfty}) and (\ref{inner3Ginfty}), using \textit{Chebfun}, matching to the constant values in (\ref{solouterGinfty}). For our particular choice of dimensionless nematic viscosities $\alpha_{2}$ and $\alpha_{3}$, all values of $\sigma_{k}^{\pm}$ (defined in (\ref{solouterGinfty})) are  close to some odd multiple of $\frac{\pi}{2}$, and thus the inner director field is largely flow-aligned and is rotated $k \pi$ times with respect to the flow direction. There are multiple choices for the outer solutions, $\sigma_{k_{1}}^{\pm}$ and $\sigma_{k_{2}}^{\pm}$, for $-1<z<0$ and $0<z<1$ respectively, yielding different asymptotic approximations.  
In Figures~\ref{fig:asymptoticsnearginft0}(a) and~\ref{fig:asymptoticsnearginft2}(a) we compare the asymptotic approximations (\ref{solouterGinfty}), (\ref{inner1Ginfty}), (\ref{inner2Ginfty}) and (\ref{inner3Ginfty}) with numerical solutions of the full system (\ref{systemstatic}) for large values of $\mathcal{G}$. The two cases are labeled as $\theta_{a_{0},\mathcal{G}}^{\ast}$ and $\theta_{\tilde{a}_{1},\mathcal{G}}^{\ast}$ respectively, depending on the initial condition used to generate them. The values of $\sigma_{k_{1}}^{\pm}$ and $\sigma_{k_{2}}^{\pm}$  are extracted from the numerical solution and used in the asymptotic approximation (\ref{solouterGinfty})-(\ref{inner3Ginfty}) (these values are different for solutions $\theta_{a_{0},\mathcal{G}}^{\ast}$ and $\theta_{\tilde{a}_{1},\mathcal{G}}^{\ast}$). Once the outer values are determined, we can compute the asymptotic approximation using the methodology outlined above. 
The asymptotic solution approximates the full numerical solution well. The asymptotic solutions also show that the boundary layers near the walls have width proportional to $\mathcal{G}^{-1/2}$, consistent with the experimental findings in Sengupta et al.\cite{Sengupta2013} In Figures \ref{fig:asymptoticsnearginft0}(b) and \ref{fig:asymptoticsnearginft2}(b), we plot the director field $\ve{n}$ associated with the equilibria $\theta_{a_{0},\mathcal{G}}^{\ast}$ and $\theta_{\tilde{a}_{1},\mathcal{G}}^{\ast}$, computed for $\mathcal{G}=100$ and $\mathcal{B}=\frac{1}{3}$. The director field is largely flow-aligned and the director field associated with $\theta_{a_{0},\mathcal{G}}^{\ast}$ exhibits a third transition layer near the centre as predicted by the asymptotic analysis.
\subsection{Equilibrium solution landscape in $\mathcal{G}$.}\label{sec:solutionlandscape}
\begin{figure}[ht!]
\centering
%
%
\begin{psfrags}%
\psfragscanon%
%
\psfrag{s03}[t][t]{\color[rgb]{0,0,0}\setlength{\tabcolsep}{0pt}\begin{tabular}{c}{\LARGE $\omega(\theta^{\ast}_{a_{n},\mathcal{G}})$}\end{tabular}}%
\psfrag{s04}[b][b]{\color[rgb]{0,0,0}\setlength{\tabcolsep}{0pt}\begin{tabular}{c}{\LARGE $\mathcal{B}$}\end{tabular}}%
\psfrag{s05}[l][l]{\color[rgb]{0,0,0}\setlength{\tabcolsep}{0pt}\begin{tabular}{l}$\mathcal{G}=$0\end{tabular}}%
\psfrag{s06}[l][l]{\color[rgb]{0,0,0}\setlength{\tabcolsep}{0pt}\begin{tabular}{l}$\mathcal{G}=$5\end{tabular}}%
\psfrag{s07}[l][l]{\color[rgb]{0,0,0}\setlength{\tabcolsep}{0pt}\begin{tabular}{l}$\mathcal{G}=$20\end{tabular}}%
\psfrag{s08}[l][l]{\color[rgb]{0,0,0}\setlength{\tabcolsep}{0pt}\begin{tabular}{l}{\LARGE Type I}\end{tabular}}%
\psfrag{s09}[l][l]{\color[rgb]{0,0,0}\setlength{\tabcolsep}{0pt}\begin{tabular}{l}{\Large increasing $\mathcal{G}$}\end{tabular}}%
%
\psfrag{x01}[t][t]{0}%
\psfrag{x02}[t][t]{0.1}%
\psfrag{x03}[t][t]{0.2}%
\psfrag{x04}[t][t]{0.3}%
\psfrag{x05}[t][t]{0.4}%
\psfrag{x06}[t][t]{0.5}%
\psfrag{x07}[t][t]{0.6}%
\psfrag{x08}[t][t]{0.7}%
\psfrag{x09}[t][t]{0.8}%
\psfrag{x10}[t][t]{0.9}%
\psfrag{x11}[t][t]{1}%
\psfrag{x12}[t][t]{-6.2832}%
\psfrag{x13}[t][t]{-4.7124}%
\psfrag{x14}[t][t]{-3.1416}%
\psfrag{x15}[t][t]{-1.5708}%
\psfrag{x16}[t][t]{0}%
\psfrag{x17}[t][t]{1.5708}%
\psfrag{x18}[t][t]{3.1416}%
\psfrag{x19}[t][t]{4.7124}%
\psfrag{x20}[t][t]{6.2832}%
%
\psfrag{v01}[r][r]{0}%
\psfrag{v02}[r][r]{0.1}%
\psfrag{v03}[r][r]{0.2}%
\psfrag{v04}[r][r]{0.3}%
\psfrag{v05}[r][r]{0.4}%
\psfrag{v06}[r][r]{0.5}%
\psfrag{v07}[r][r]{0.6}%
\psfrag{v08}[r][r]{0.7}%
\psfrag{v09}[r][r]{0.8}%
\psfrag{v10}[r][r]{0.9}%
\psfrag{v11}[r][r]{1}%
\psfrag{v12}[r][r]{0}%
\psfrag{v13}[r][r]{0.1}%
\psfrag{v14}[r][r]{0.2}%
\psfrag{v15}[r][r]{0.3}%
\psfrag{v16}[r][r]{0.4}%
\psfrag{v17}[r][r]{0.5}%
\psfrag{v18}[r][r]{0.6}%
\psfrag{v19}[r][r]{0.7}%
\psfrag{v20}[r][r]{0.8}%
\psfrag{v21}[r][r]{0.9}%
\psfrag{v22}[r][r]{1}%
%
\resizebox{8cm}{!}{\includegraphics{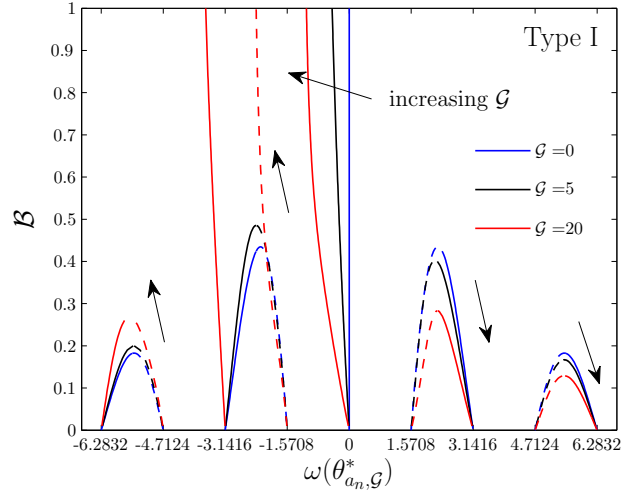}}%
\end{psfrags}%
%

\caption{Evolution of the steady-state solutions of Type I as $\mathcal{G}$ increases. The solid and dashed lines represent, respectively, the values of $\omega(\theta_{a_{n},\mathcal{G}}^{\ast})$ for which the steady states, $\theta_{a_{n},\mathcal{G}}^{\ast}$, are stable or unstable.}
\label{fig:solutionlandscapetype1}
\end{figure}

\begin{figure}[ht!]
\centering
%
%
\begin{psfrags}%
\psfragscanon%
%
\psfrag{s03}[t][t]{\color[rgb]{0,0,0}\setlength{\tabcolsep}{0pt}\begin{tabular}{c}{\LARGE $\omega(\theta_{\tilde{a}_{n},\mathcal{G}}^{\ast})$}\end{tabular}}%
\psfrag{s04}[b][b]{\color[rgb]{0,0,0}\setlength{\tabcolsep}{0pt}\begin{tabular}{c}{\LARGE $\mathcal{B}$}\end{tabular}}%
\psfrag{s05}[l][l]{\color[rgb]{0,0,0}\setlength{\tabcolsep}{0pt}\begin{tabular}{l}$\mathcal{G}=$0\end{tabular}}%
\psfrag{s06}[l][l]{\color[rgb]{0,0,0}\setlength{\tabcolsep}{0pt}\begin{tabular}{l}$\mathcal{G}=$5\end{tabular}}%
\psfrag{s07}[l][l]{\color[rgb]{0,0,0}\setlength{\tabcolsep}{0pt}\begin{tabular}{l}$\mathcal{G}=$20\end{tabular}}%
\psfrag{s08}[l][l]{\color[rgb]{0,0,0}\setlength{\tabcolsep}{0pt}\begin{tabular}{l}{\LARGE Type II}\end{tabular}}%
\psfrag{s09}[l][l]{\color[rgb]{0,0,0}\setlength{\tabcolsep}{0pt}\begin{tabular}{l}{\Large increasing $\mathcal{G}$}\end{tabular}}%
%
\psfrag{x01}[t][t]{0}%
\psfrag{x02}[t][t]{0.1}%
\psfrag{x03}[t][t]{0.2}%
\psfrag{x04}[t][t]{0.3}%
\psfrag{x05}[t][t]{0.4}%
\psfrag{x06}[t][t]{0.5}%
\psfrag{x07}[t][t]{0.6}%
\psfrag{x08}[t][t]{0.7}%
\psfrag{x09}[t][t]{0.8}%
\psfrag{x10}[t][t]{0.9}%
\psfrag{x11}[t][t]{1}%
\psfrag{x12}[t][t]{-2}%
\psfrag{x13}[t][t]{-1.5}%
\psfrag{x14}[t][t]{-1}%
\psfrag{x15}[t][t]{-0.5}%
\psfrag{x16}[t][t]{0}%
\psfrag{x17}[t][t]{0.5}%
\psfrag{x18}[t][t]{1}%
\psfrag{x19}[t][t]{1.5}%
\psfrag{x20}[t][t]{2}%
%
\psfrag{v01}[r][r]{0}%
\psfrag{v02}[r][r]{0.1}%
\psfrag{v03}[r][r]{0.2}%
\psfrag{v04}[r][r]{0.3}%
\psfrag{v05}[r][r]{0.4}%
\psfrag{v06}[r][r]{0.5}%
\psfrag{v07}[r][r]{0.6}%
\psfrag{v08}[r][r]{0.7}%
\psfrag{v09}[r][r]{0.8}%
\psfrag{v10}[r][r]{0.9}%
\psfrag{v11}[r][r]{1}%
\psfrag{v12}[r][r]{0}%
\psfrag{v13}[r][r]{0.5}%
\psfrag{v14}[r][r]{1}%
\psfrag{v15}[r][r]{1.5}%
\psfrag{v16}[r][r]{2}%
%
\resizebox{8cm}{!}{\includegraphics{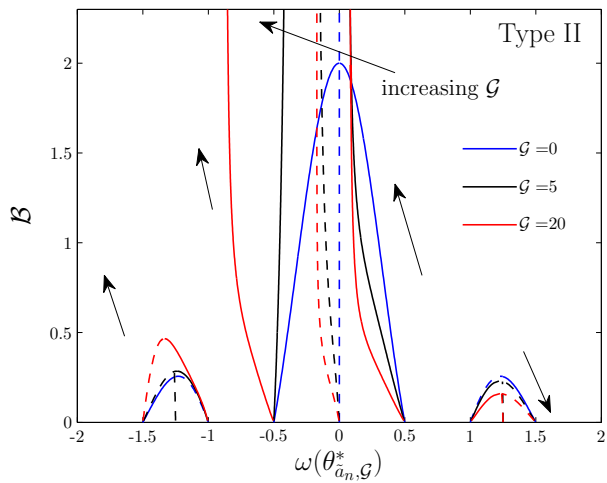}}%
\end{psfrags}%
%

\caption{Evolution of the steady-state solutions of Type II as $\mathcal{G}$ increases. The solid and dashed lines represent, respectively, the values of $\omega(\theta_{\tilde{a}_{n},\mathcal{G}}^{\ast})$ for which the steady states, $\theta_{\tilde{a}_{n},\mathcal{G}}^{\ast}$, are stable or unstable.}
\label{fig:solutionlandscapetype2}
\end{figure}

In this section, we study how the static solution landscape for the system (\ref{systemtdep}) evolves as the pressure gradient $\mathcal{G}$ increases. In Section~\ref{sec:G0}, we compute the static equilibria, $\theta^*_{a}$ for $\mathcal{G}=0$. In what follows, we let $\theta^*_{a,\mathcal{G}}$ denote the numerically computed equilibrium, via continuation methods with $\theta_{a}^{\ast}$ as initial condition.
We numerically compute the stability of the equilibria with $\mathcal{G}>0$ (using the function \textit{eigs} of the MATLAB package \textit{Chebfun}) and find that the stability properties of the $\mathcal{G}=0$ equilibria propagate to the $\mathcal{G}>0$ cases. Figures \ref{fig:solutionlandscapetype1} and \ref{fig:solutionlandscapetype2} show the evolution of the steady state solutions, $\theta^{\ast}_{a_{n}}$ and $\theta^{\ast}_{\tilde{a}_{n}}$, as $\mathcal{G}$ increases. 
For $\mathcal{G}=0$ and $\mathcal{B}>\mathcal{B}^{\ast}_{1}$, the trivial solution $\theta_{a_{0}}^{\ast} \equiv 0$ is the unique stable equilibrium. For $\mathcal{G}>0$ the trivial solution is not an equilibrium and for $\mathcal{B} > \mathcal{B}^{\ast}_{1}$, $\theta_{a_{0},\mathcal{G}}^{\ast}$ is not the unique stable equilibrium. As the pressure gradient $\mathcal{G}$ increases, new equilibria appear for $\mathcal{B} > \mathcal{B}_{1}^{\ast}$. Additionally, some equilibria, e.g. those with a large positive winding number, become suppressed or have a smaller window of existence in $\mathcal{B}$, as $\mathcal{G}$ increases.

We believe that the asymmetry in the solution branches with positive and negative winding numbers for $\mathcal{G}>0$ is a consequence of the fact that we work with unit-vector fields, and not director fields without a direction. We speculate that a more sophisticated model, such as the Beris--Edwards model for nematodynamics which accounts for the head--tail symmetry of nematic molecules, may resolve this asymmetry between positive and negative winding numbers for large $\mathcal{G}$. 

\begin{figure}[ht!]
\centering
%
%
\begin{psfrags}%
\psfragscanon%
%
\psfrag{s03}[t][t]{\color[rgb]{0,0,0}\setlength{\tabcolsep}{0pt}\begin{tabular}{c}{\LARGE $\mathcal{G}$}\end{tabular}}%
\psfrag{s04}[b][b]{\color[rgb]{0,0,0}\setlength{\tabcolsep}{0pt}\begin{tabular}{c}{\LARGE $\mathcal{B}^{\ast}_{i,\mathcal{G}}$}\end{tabular}}%
\psfrag{s05}[l][l]{\color[rgb]{0,0,0}\setlength{\tabcolsep}{0pt}\begin{tabular}{l}{\Large $\mathcal{B}_{-2,\mathcal{G}}^{\ast}$}\end{tabular}}%
\psfrag{s06}[l][l]{\color[rgb]{0,0,0}\setlength{\tabcolsep}{0pt}\begin{tabular}{l}{\Large $\mathcal{B}_{-3,\mathcal{G}}^{\ast}$}\end{tabular}}%
\psfrag{s07}[l][l]{\color[rgb]{0,0,0}\setlength{\tabcolsep}{0pt}\begin{tabular}{l}{\Large $\mathcal{B}_{-4,\mathcal{G}}^{\ast}$}\end{tabular}}%
\psfrag{s08}[l][l]{\color[rgb]{0,0,0}\setlength{\tabcolsep}{0pt}\begin{tabular}{l}{\Large $\mathcal{B}_{2,\mathcal{G}}^{\ast}$}\end{tabular}}%
\psfrag{s09}[l][l]{\color[rgb]{0,0,0}\setlength{\tabcolsep}{0pt}\begin{tabular}{l}{\Large $\mathcal{B}_{4,\mathcal{G}}^{\ast}$}\end{tabular}}%
\psfrag{s10}[l][l]{\color[rgb]{0,0,0}\setlength{\tabcolsep}{0pt}\begin{tabular}{l}{\Large $\mathcal{B}_{3,\mathcal{G}}^{\ast}$}\end{tabular}}%
%
\psfrag{x01}[t][t]{0}%
\psfrag{x02}[t][t]{0.1}%
\psfrag{x03}[t][t]{0.2}%
\psfrag{x04}[t][t]{0.3}%
\psfrag{x05}[t][t]{0.4}%
\psfrag{x06}[t][t]{0.5}%
\psfrag{x07}[t][t]{0.6}%
\psfrag{x08}[t][t]{0.7}%
\psfrag{x09}[t][t]{0.8}%
\psfrag{x10}[t][t]{0.9}%
\psfrag{x11}[t][t]{1}%
\psfrag{x12}[t][t]{0}%
\psfrag{x13}[t][t]{5}%
\psfrag{x14}[t][t]{10}%
\psfrag{x15}[t][t]{15}%
\psfrag{x16}[t][t]{20}%
\psfrag{x17}[t][t]{25}%
\psfrag{x18}[t][t]{30}%
\psfrag{x19}[t][t]{35}%
%
\psfrag{v01}[r][r]{0}%
\psfrag{v02}[r][r]{0.1}%
\psfrag{v03}[r][r]{0.2}%
\psfrag{v04}[r][r]{0.3}%
\psfrag{v05}[r][r]{0.4}%
\psfrag{v06}[r][r]{0.5}%
\psfrag{v07}[r][r]{0.6}%
\psfrag{v08}[r][r]{0.7}%
\psfrag{v09}[r][r]{0.8}%
\psfrag{v10}[r][r]{0.9}%
\psfrag{v11}[r][r]{1}%
\psfrag{v12}[r][r]{0}%
\psfrag{v13}[r][r]{0.5}%
\psfrag{v14}[r][r]{1}%
\psfrag{v15}[r][r]{1.5}%
\psfrag{v16}[r][r]{2}%
\psfrag{v17}[r][r]{2.5}%
%
\resizebox{8cm}{!}{\includegraphics{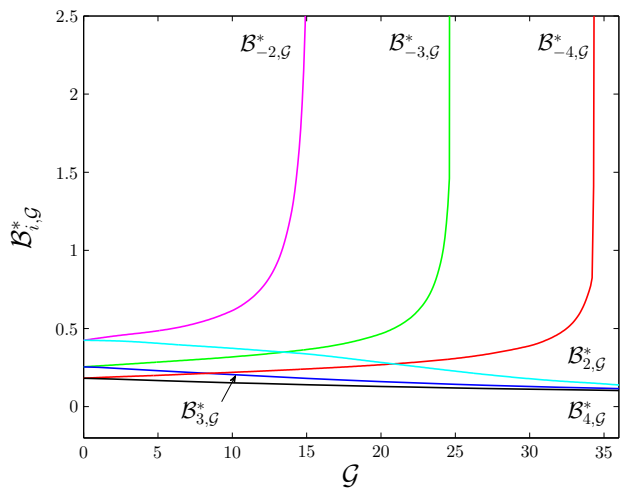}}%
\end{psfrags}%
%

\caption{Evolution of the critical values $\mathcal{B}_{i,\mathcal{G}}^{\ast}$ as $\mathcal{G}$ increases.}
\label{fig:bast}
\end{figure}

Let $\mathcal{B}_{i, \mathcal{G}}^{\ast}$ denote a critical value of $\mathcal{B}$ for a fixed $\mathcal{G} > 0$; this definition is analogous to the definition of $\mathcal{B}^{\ast}_{i}$ for $\mathcal{G} = 0$. We conjecture that there is a saddle-node bifurcation at each critical value such that if $n>0$, the stable branch, $\theta^{\ast}_{a_{2n},\mathcal{G}}$, and the unstable branch, $\theta^{\ast}_{a_{2n-1},\mathcal{G}}$ ($\theta^{\ast}_{a_{2n+1},\mathcal{G}}$ for $n<0$), collide at $\mathcal{B}=\mathcal{B}^{\ast}_{2n,\mathcal{G} }$ and cease to exist for $\mathcal{B}> \mathcal{B}^{\ast}_{2n,\mathcal{G}}$ (similarly for $\mathcal{B}^{\ast}_{2n+1, \mathcal{G}}$ and solutions of Type II). In Figure~\ref{fig:bast} we plot the critical values $\mathcal{B}^{\ast}_{i,\mathcal{G}}$ $i=\pm 2, 3, \ldots$ as a function of the pressure gradient. For example, if $\mathcal{G} \approx 15$, the critical value $\mathcal{B}_{-2,\mathcal{G}}^{\ast} \rightarrow \infty$ so that for $\mathcal{G}>15$, the solution branches $\theta^{\ast}_{a_{-2},\mathcal{G}}$ and $\theta^{\ast}_{a_{-1},\mathcal{G}}$ do not coalesce and exist for all $\mathcal{B}$. 
\section{Time-dependent solutions}\label{sec:time}
\begin{figure}[ht!]
\centering
%
%
\begin{psfrags}%
\psfragscanon%
%
\psfrag{s03}[t][t]{\color[rgb]{0,0,0}\setlength{\tabcolsep}{0pt}\begin{tabular}{c}{\LARGE $C$}\end{tabular}}%
\psfrag{s04}[b][b]{\color[rgb]{0,0,0}\setlength{\tabcolsep}{0pt}\begin{tabular}{c}{\LARGE $\omega(\theta^{\ast})$}\end{tabular}}%
\psfrag{s05}[l][l]{\color[rgb]{0,0,0}\setlength{\tabcolsep}{0pt}\begin{tabular}{l}$\omega(\theta^{\ast}_{a_{-6},2})$\end{tabular}}%
\psfrag{s06}[l][l]{\color[rgb]{0,0,0}\setlength{\tabcolsep}{0pt}\begin{tabular}{l}$\omega(\theta^{\ast}_{a_{-4},2})$\end{tabular}}%
\psfrag{s07}[l][l]{\color[rgb]{0,0,0}\setlength{\tabcolsep}{0pt}\begin{tabular}{l}$\omega(\theta^{\ast}_{a_{-2},2})$\end{tabular}}%
\psfrag{s08}[l][l]{\color[rgb]{0,0,0}\setlength{\tabcolsep}{0pt}\begin{tabular}{l}$\omega(\theta^{\ast}_{a_{0},2})$\end{tabular}}%
\psfrag{s09}[l][l]{\color[rgb]{0,0,0}\setlength{\tabcolsep}{0pt}\begin{tabular}{l}$\omega(\theta^{\ast}_{a_{2},2})$\end{tabular}}%
\psfrag{s10}[l][l]{\color[rgb]{0,0,0}\setlength{\tabcolsep}{0pt}\begin{tabular}{l}$\omega(\theta^{\ast}_{a_{4},2})$\end{tabular}}%
\psfrag{s11}[l][l]{\color[rgb]{0,0,0}\setlength{\tabcolsep}{0pt}\begin{tabular}{l}$\omega(\theta^{\ast}_{a_{6},2})$\end{tabular}}%
%

\psfrag{x01}[t][t]{-$3\pi$}%
\psfrag{x02}[t][t]{-$2\pi$}%
\psfrag{x03}[t][t]{-$\pi$}%
\psfrag{x04}[t][t]{{\Large $\textcolor{red}{C^{\ast}}$}}%
\psfrag{x05}[t][t]{0}%
\psfrag{x06}[t][t]{$\pi$}%
\psfrag{x07}[t][t]{$2\pi$}%
\psfrag{x08}[t][t]{$3\pi$}%
%
\psfrag{v01}[r][r]{-3}%
\psfrag{v02}[r][r]{-2}%
\psfrag{v03}[r][r]{-1}%
\psfrag{v04}[r][r]{0}%
\psfrag{v05}[r][r]{1}%
\psfrag{v06}[r][r]{2}%
\psfrag{v07}[r][r]{3}%
%
\resizebox{8cm}{!}{\includegraphics{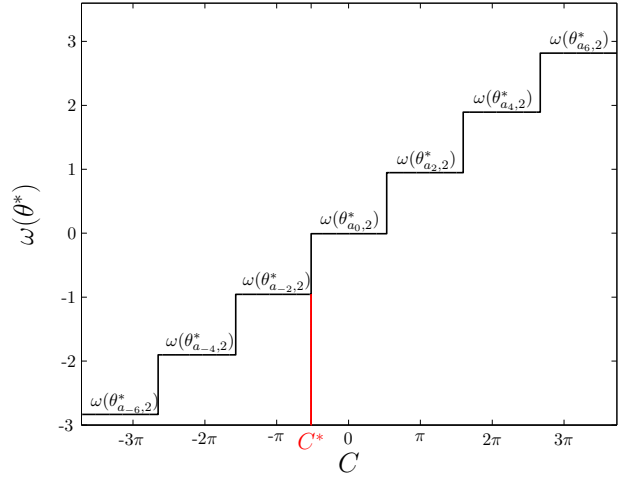}}%
\end{psfrags}%
%

\caption{Winding number for the solution of the system (\ref{systemtdep}) with $\mathcal{B}=\frac{1}{10}$, $\mathcal{G}=2$, with different linear initial conditions $\Theta(z)=Cz$, $C \in [-7\frac{\pi}{2},7\frac{\pi}{2}]$. The critical value $C^{\ast}$ is indicated on the $x$-axis.}
\label{fig:winding}
\end{figure}

In this section, we study the time-dependent behavior of the system (\ref{systemtdep}). We numerically compute the time-dependent solutions using a self-implemented finite-difference method, with mesh resolution $\Delta z= $0.0125 and time step $\Delta t=0.01$. As we have seen in Section \ref{sec:equilibrium}, there are multiple static equilibria for a given pair $\left(\mathcal{G}, \mathcal{B} \right)$ and it is of interest to investigate steady-state selection, for different choices of the initial conditions.
We perform a preliminary investigation of the parameter space by working with either constant or linear initial conditions. We conclude that the time-dependent system converges to:
\begin{align}
\theta^{\ast}_{a_{0}, \mathcal{G}} & \mbox{ if } & \Theta(z)= C, \label{initial:constant} \\
\theta^{\ast}_{a_{n},\mathcal{G}} & \mbox{ if } & \Theta(z)=Cz,\label{initial:lineart1} \\
\theta^{\ast}_{\tilde{a}_{n},\mathcal{G}} & \mbox{ if } & \Theta(z)=Cz+\frac{\pi}{2}, \label{initial:lineart2}
\end{align}
where $C$ is a constant. We note that the initial conditions in \eqref{initial:constant}--\eqref{initial:lineart2} do not satisfy the boundary conditions in (\ref{systemtdep}) and in Section \ref{sec:tuning}, we propose alternative initial conditions that respect these boundary conditions.
In Figure \ref{fig:winding} we use linear initial conditions \eqref{initial:lineart1} that have $C \in [-\frac{7\pi}{2},\frac{7\pi}{2}]$, $\mathcal{G}=2$, $\mathcal{B}=\frac{1}{10}$, and find that the steady state converges to different equilibria $\theta^{\ast}_{a_{n},2}$, depending on the initial value $C$. We compute the corresponding winding numbers and use the winding number to label the static equilibria in Figure \ref{fig:winding}. 
Particularly, for any pair ($\mathcal{G}$,$\mathcal{B}$), we numerically find a critical value $C^{\ast}$ such that if $C \in (C^{\ast}-\epsilon,C^{\ast}+\epsilon)$, with $\epsilon>0$ sufficiently small, we have
\begin{equation}\label{Cast}
\lim_{t \rightarrow \infty} \theta(t,z;Cz)= \left\{\begin{array}{l l} \theta_{a_{-2},\mathcal{G}}^{\ast} & \text{ if } C \in (C^{\ast}-\epsilon,C^{\ast}), \\ \theta_{a_{0},\mathcal{G}}^{\ast} & \text{ if } C \in [C^{\ast},C^{\ast}+\epsilon). \end{array}\right.
\end{equation}\\
\begin{figure}[ht!]
\centering
%
%
\begin{psfrags}%
\psfragscanon%
%
\psfrag{s04}[t][t]{\color[rgb]{0,0,0}\setlength{\tabcolsep}{0pt}\begin{tabular}{c}{\LARGE $z$}\end{tabular}}%
\psfrag{s05}[l][l]{\color[rgb]{0,0,0}\setlength{\tabcolsep}{0pt}\begin{tabular}{l}{\Large $\theta^{\ast}_{a_{0},\mathcal{G}}$}\end{tabular}}%
\psfrag{s06}[l][l]{\color[rgb]{0,0,0}\setlength{\tabcolsep}{0pt}\begin{tabular}{l}{\Large $\theta^{\ast}_{a_{-2},\mathcal{G}}$}\end{tabular}}%
\psfrag{s07}[l][l]{\color[rgb]{0,0,0}\setlength{\tabcolsep}{0pt}\begin{tabular}{l}{\Large $\Theta= C^{\ast}z$}\end{tabular}}%
\psfrag{s08}[l][l]{\color[rgb]{1,0,0}\setlength{\tabcolsep}{0pt}\begin{tabular}{l}{\Large System \eqref{systemtdep} approaches $\theta_{a_{-2},\mathcal{G}}^{\ast}$}\end{tabular}}%
\psfrag{s09}[l][l]{\color[rgb]{0,0,1}\setlength{\tabcolsep}{0pt}\begin{tabular}{l}{\Large System \eqref{systemtdep} approaches $\theta_{a_{0},\mathcal{G}}^{\ast}$}\end{tabular}}%
%
\psfrag{x01}[t][t]{0}%
\psfrag{x02}[t][t]{0.1}%
\psfrag{x03}[t][t]{0.2}%
\psfrag{x04}[t][t]{0.3}%
\psfrag{x05}[t][t]{0.4}%
\psfrag{x06}[t][t]{0.5}%
\psfrag{x07}[t][t]{0.6}%
\psfrag{x08}[t][t]{0.7}%
\psfrag{x09}[t][t]{0.8}%
\psfrag{x10}[t][t]{0.9}%
\psfrag{x11}[t][t]{1}%
\psfrag{x12}[t][t]{-1}%
\psfrag{x13}[t][t]{-0.8}%
\psfrag{x14}[t][t]{-0.6}%
\psfrag{x15}[t][t]{-0.4}%
\psfrag{x16}[t][t]{-0.2}%
\psfrag{x17}[t][t]{0}%
\psfrag{x18}[t][t]{0.2}%
\psfrag{x19}[t][t]{0.4}%
\psfrag{x20}[t][t]{0.6}%
\psfrag{x21}[t][t]{0.8}%
\psfrag{x22}[t][t]{1}%
%
\psfrag{v01}[r][r]{0}%
\psfrag{v02}[r][r]{0.1}%
\psfrag{v03}[r][r]{0.2}%
\psfrag{v04}[r][r]{0.3}%
\psfrag{v05}[r][r]{0.4}%
\psfrag{v06}[r][r]{0.5}%
\psfrag{v07}[r][r]{0.6}%
\psfrag{v08}[r][r]{0.7}%
\psfrag{v09}[r][r]{0.8}%
\psfrag{v10}[r][r]{0.9}%
\psfrag{v11}[r][r]{1}%
\psfrag{v12}[r][r]{-3}%
\psfrag{v13}[r][r]{-2}%
\psfrag{v14}[r][r]{-1}%
\psfrag{v15}[r][r]{0}%
\psfrag{v16}[r][r]{1}%
\psfrag{v17}[r][r]{2}%
\psfrag{v18}[r][r]{3}%
%
\resizebox{8cm}{!}{\includegraphics{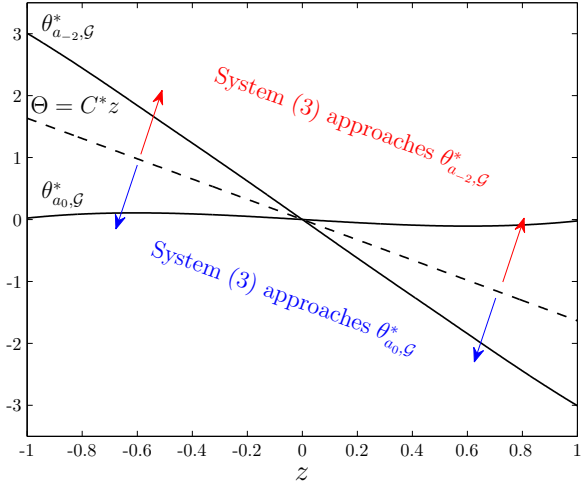}}%
\end{psfrags}%
%

\caption{Solutions $\theta^{\ast}_{a_{0},\mathcal{G}}$ and $\theta^{\ast}_{a_{2},\mathcal{G}}$ obtained with $\mathcal{B}=\frac{1}{10}$ and $\mathcal{G}=2$. The critical initial condition $\Theta(z)=C^{\ast}z$ is plotted with dashed line.}
\label{fig:Cast}
\end{figure}
Figure \ref{fig:Cast} plots the initial condition $\Theta(z)=C^{\ast}(z)$, where $C^{\ast}$ is the critical value obtained with $\mathcal{G}=2$ and $\mathcal{B}= \frac{1}{10}$. System (\ref{systemtdep}) with initial condition $\Theta(z)= C z$ approaches either $\theta^{\ast}_{a_{0},2}$ or $\theta^{\ast}_{a_{-2},2}$ if $C\geq C^{\ast}$ or $C< C^{\ast}$, respectively.   
\subsection{Tuning the pressure gradient and the boundary conditions}\label{sec:tuning}
The pressure gradient and boundary conditions have been assumed to be constants in our computations to this point. However, it is of experimental interest to consider situations where both the pressure gradient and boundary conditions are continuously tuned over a short period of time until they attain the desired state. We consider tuning the flow at a rate $\delta$ by applying 
\begin{equation}\label{tuneG}
\mathcal{G}(t) = \left \{ \begin{array}{l c l}
0 & & \mbox{ if } t\leq t_{1},
\\
\bar{\mathcal{G}}\tanh(\delta (t-t_{1})) & & \mbox{otherwise}.
\end{array}\right.
\end{equation}
Similarly, we apply time-dependent anchoring conditions of the form 
\begin{equation}\label{tuneBound}
\begin{array}{r l}
\theta_{z}(1,t)= & \left \{ \begin{array}{l c l} C & &\mbox{ if } t\leq t_{2}  \\ \\ C(1-\tanh(\kappa (t-t_{2}))) & & \\
 \hspace*{1cm}-\displaystyle\frac{\sin(2\theta(1,t))\tanh(\kappa(t-t_{2}))}{\mathcal{B}} & &\mbox{otherwise},
\end{array}\right.
\\
\\
\theta_{z}(-1,t)= &  \left\{ \begin{array}{l c l}
 C  & & \mbox{ if } t\leq t_{2} \\ \\
C(1-\tanh(\kappa (t-t_{2}))) & & \\ \hspace*{1cm} +\displaystyle \frac{\sin(2\theta(1,t))\tanh(\kappa(t-t_{2}))}{\mathcal{B}} & & \mbox{otherwise},
\end{array}\right.
\end{array}
\end{equation}
for some constant $\kappa>0$. In particular, these conditions are satisfied by the initial (linear) condition $\Theta=Cz$ for $t\leq t_2$ and then, the anchoring is switched on with a tuning rate $\kappa$, to attain the required weak anchoring conditions at $z=\pm 1$.

We numerically study this modified dynamic system, using \eqref{tuneG} and \eqref{tuneBound}, and find that if $t_{1} \leq t_{2}$, then the final steady state is identical to the steady state attained with constant values $\mathcal{G}=\bar{\mathcal{G}}$ and boundary conditions \eqref{boundtop}--\eqref{boundbot}, for the parameter sweep that we performed. This indicates that if we first apply a pressure gradient and then enforce strong anchoring, the system will always relax to the same equilibrium state, regardless of the time delay between application of the pressure gradient and anchoring.   

On the other hand, if we apply the anchoring condition before the pressure gradient by choosing $t_1>t_2$, then a different steady state can be attained, depending on the time delay and the respective rates. As an illustrative example, we find that if $\Theta=Cz$ with $C < C^{\ast}$ and $\mathcal{B}> \mathcal{B}_{-2}^{\ast}$, solutions of system (\ref{systemtdep}) with (\ref{tuneG})--(\ref{tuneBound}) may approach the equilibrium solution, $\theta_{a_{0}, \bar{\mathcal{G}}}^{\ast}$, instead of the expected solution, $\theta_{a_{-2},\bar{\mathcal{G}}}^{\ast}$. This can be explained as follows: when $t_{\rm 2}< t\leq t_{\rm 1}$, i.e.~while $\mathcal{G}=0$, the trivial solution $\theta_{a_{0}}^{\ast} = 0$ is the unique steady state and thus the system must approach this solution during the early stages. As a consequence, when the flow begins ($t> t_{\rm 1}$), the solution is already sufficiently close to $\theta_{a_{0}}^{\ast}$ and thus can no longer access the equilibrium state $\theta_{a_{-2},\bar{\mathcal{G}}}^{\ast}$, as it would do if $t_{\rm 1}\leq t_{\rm 2}$. Hence, given model parameters $\bar{\mathcal{G}}$, $\mathcal{B}$, $t_{\rm 2}$, $\kappa$ and $\delta$, if the initial condition is $\Theta= Cz$, one can define a critical value $t_{\rm 1}^{\ast}(C)$ such that
\begin{equation}\label{t2ast}
\left \{ \begin{array}{l r} \displaystyle\lim_{t \rightarrow \infty} \theta(t,z;Cz)= \theta_{a_{-2},\bar{\mathcal{G}}}^{\ast} & \text{ if } t_{\rm 1}<t_{\rm 1}^{\ast}\\
\\
\displaystyle\lim_{t \rightarrow \infty} \theta(t,z;Cz)= \theta_{a_{0},\bar{\mathcal{G}}}^{\ast}  & \text{ if } t_{\rm 1}\geq t_{\rm 1}^{\ast}. \end{array} \right \}
\end{equation}
If $C$ is such that $\lim_{t \rightarrow \infty} \theta(t,z;Cz)= \theta_{a_{-2},\bar{\mathcal{G}}}^{\ast}$ for all $t_{1}>0$, $t_{1}^{\ast}(C)$ is not defined.
\begin{figure}[ht!]
\centering
%
%
\begin{psfrags}%
\psfragscanon%
%
\psfrag{s03}[t][t]{\color[rgb]{0,0,0}\setlength{\tabcolsep}{0pt}\begin{tabular}{c}{\LARGE $C$}\end{tabular}}%
\psfrag{s04}[b][b]{\color[rgb]{0,0,0}\setlength{\tabcolsep}{0pt}\begin{tabular}{c}{\LARGE $t_{1}^{\ast}(C)$}\end{tabular}}%
\psfrag{s05}[l][l]{\color[rgb]{0,0,0}\setlength{\tabcolsep}{0pt}\begin{tabular}{l}{\Large $\mathcal{B}=0.5$}\end{tabular}}%
\psfrag{s06}[l][l]{\color[rgb]{0,0,0}\setlength{\tabcolsep}{0pt}\begin{tabular}{l}{\Large $\mathcal{B}=0.8$}\end{tabular}}%
\psfrag{s07}[l][l]{\color[rgb]{0,0,0}\setlength{\tabcolsep}{0pt}\begin{tabular}{l}{\Large $\mathcal{B}=1$}\end{tabular}}%
%
\psfrag{x01}[t][t]{-4}%
\psfrag{x02}[t][t]{-3.8}%
\psfrag{x03}[t][t]{-3.6}%
\psfrag{x04}[t][t]{-3.4}%
\psfrag{x05}[t][t]{-3.2}%
\psfrag{x06}[t][t]{-3}%
\psfrag{x07}[t][t]{-2.8}%
\psfrag{x08}[t][t]{-2.6}%
\psfrag{x09}[t][t]{-2.4}%
\psfrag{x10}[t][t]{-2.2}%
\psfrag{x11}[t][t]{-2}%
%
\psfrag{v01}[r][r]{0}%
\psfrag{v02}[r][r]{0.1}%
\psfrag{v03}[r][r]{0.2}%
\psfrag{v04}[r][r]{0.3}%
\psfrag{v05}[r][r]{0.4}%
\psfrag{v06}[r][r]{0.5}%
\psfrag{v07}[r][r]{0.6}%
\psfrag{v08}[r][r]{0.7}%
\psfrag{v09}[r][r]{0.8}%
%
\resizebox{8cm}{!}{\includegraphics{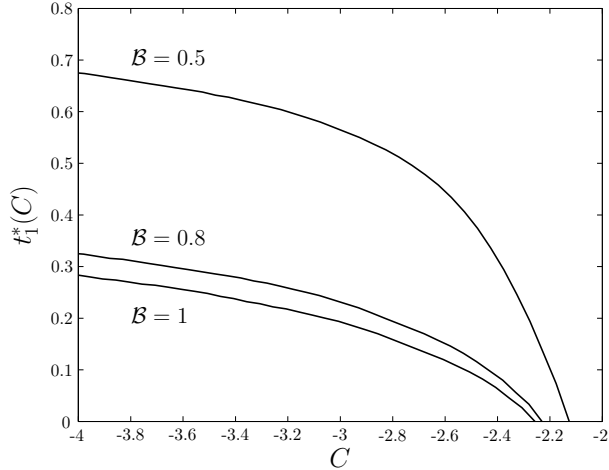}}%
\end{psfrags}%
%

\caption{Critical values $t_{\rm 1}^{\ast}(C)$ obtained when solving the system (\ref{systemtdep}),(\ref{tuneG})--(\ref{tuneBound}) with $\bar{\mathcal{G}}=40$, $t_{\rm 2}=0$ and $\delta=\kappa =5$. For $t_1 < t_1^{\ast}(C)$, the solution evolves to the steady state $\theta_{a_{-2},\bar{\mathcal{G}}}^{\ast}$; for $t_1 \geq t_1^{\ast}(C)$, the system evolves to the steady state $\theta_{a_{0},\bar{\mathcal{G}}}^{\ast}$.  Note that $t_{1}^{\ast}=0$ when $C=C^{\ast}$ (see Definition \eqref{Cast}).}
\label{fig:howlateIcanstart}
\end{figure}

Figure \ref{fig:howlateIcanstart} shows the dependence of the critical times $t_{\rm 1}^{\ast}$ on $C$ and $\mathcal{B}$.  We observe that, as the inverse anchoring strength $\mathcal{B}$ increases, the critical time $t_{1}^{\ast}\left(C\right)$ decreases. 
This is expected since as $\mathcal{B}$ increases, the anchoring strength decreases and thus the system is able to reorient itself more easily.

\section{Conclusions}
\label{sec:conclusions}
We have explored the static equilibria of a Leslie--Ericksen model for a unidirectional uniaxial nematic flow in a prototype microfluidic channel, as a function of the pressure gradient $\mathcal{G}$ and inverse anchoring strength, $\mathcal{B}$. As $\mathcal{B} \to 0$, we approach the strong-anchoring limit. In particular, the weak- and strong-flow solutions (obtained with weak anchoring) in Anderson et al.\cite{Linda2015} correspond to $\theta^{\ast}_{a_{0},\mathcal{G}}$ and $\theta_{\tilde{a}_{1},\mathcal{G}}$. As $\mathcal{B}\to 0$, the solution $\theta_{\tilde{a}_{1},\mathcal{G}}$ has $1/2$-winding number (associated with a rotation by $\pi$ radians between the top and bottom surfaces) consistent with the Dirichlet conditions for the strong-flow solution in Anderson et al.\cite{Linda2015} Our stability analysis suggests that both solutions are stable when $\mathcal{G}=0$ and do not lose stability as $\mathcal{G}$ increases. In Anderson et al.\cite{Linda2015} the authors report that the strong-flow solution has lower energy than the weak-flow solution for large $\mathcal{G}$ and the critical $\mathcal{G}^*$ depends on the anchoring strength. This is in line with our stability analysis and we speculate that the unstable solution branches in our numerical picture may provide valuable information about how the different solution branches are connected in the full solution landscape.

We assume symmetry in the flow profile, which allows the liquid crystal dynamics to be decouple from the flow dynamics. This enables us to determine explicit and asymptotic solutions that provide key insight into the system behavior. In practice we might expect to observe transitions between the steady states that we have computed here. However, in evolving from one steady state to another the molecules may assume configurations that do not exhibit symmetry, and so this behavior will not be captured by our model. Solving the fully coupled Leslie--Ericksen model would determine the range of validity of our model in such situations.

We numerically find static equilibria with large winding numbers that are linearly stable within the simple Leslie--Ericksen model. We expect these distorted equilibria to lose stability with respect to perturbations in the $x$ and $y$ directions and so are unlikely to be observable in practice. Finally, we perform a preliminary investigation of the sensitivity of dynamic solutions to initial conditions. Working with a linear initial condition, we numerically find critical values that separate basins of attraction for the distinct steady states. Further, we also study the effect of varying the pressure gradient and anchoring conditions with time and how the rate of change can affect the critical initial conditions that lead to the selection of a particular steady state. This numerical experiment may guide future physical experiments on these lines if experimentalists can control fluid flow and anchoring conditions with time, so as to attain a desired state or at least control transient dynamics. 
We hope that our results may aid experimentalists to design new control strategies for microfluidic transport and mixing phenomena. 
\section*{Acknowledgments}
This work was carried out thanks to the financial support of the ``Spanish Ministry of Economy and Competitiveness'' under projects MTM2011-22658 and MTM2015-64865-P. The authors gratefully acknowledge many helpful discussions with L.~J.~Cummings and D.~Vigolo, and discussions at an Oxford Collaborative Workshop Initiative workshop. IMG gratefully acknowledges support from the Royal Society through a University Research Fellowship. AM's research is 
supported by an EPSRC Career Acceleration Fellowship EP/J00\-1686/1 and EP/J001686/2, an OCIAM
Visiting Fellowship and the Advanced Studies Centre at Keble College. AM and IMG are grateful for discussions on nematic microfluidics with A. Sengupta.

\newpage
\appendix
\section{Leslie--Ericksen continuum theory for nematodynamics}\label{ericksenmodel}
The Leslie--Ericksen dynamic theory is widely accepted to model dynamic phenomena in nematic liquid crystals. A unit vector $\ve{n}=(n_{1},n_{2},n_{3})$, called the \textit{director}, is defined to describe the local direction of the average molecular alignment in liquid crystals, while the instantaneous motion of the fluid is described by its velocity vector $\ve{v}=(v_{1},v_{2},v_{3})$. The full equations for nemato-dynamics describe the evolution of $\ve{n}$ and $\ve{v}$. When electromagnetic and gravitational forces are disregarded, the Leslie--Ericksen model for incompressible fluids is:\cite{Leslie1968,LeslieExistence,LeslieGeneral}
\begin{subequations}\label{le}
\begin{alignat}{2}
v_{i,i}&=0 && \mbox{ in }\Omega, \label{le1}
\\
\rho \frac{{\rm d} v_{i}}{{\rm d} \hat{t}}&= \sigma_{ji,j} &&\mbox{ in } \Omega \times (0,+\infty), \label{le2}
\\
\rho_{1} \frac{{\rm d}}{{\rm d} \hat{t}}\big(\frac{{\rm d} n_{i}}{{\rm d} \hat{t}}+ \ve{v}\cdot \nabla n_{i} \big)&= g_{i} + \pi_{ji,j} &&\mbox{ in } \Omega \times (0,+\infty), \label{le3}
\end{alignat}
\end{subequations}
where $\xi_{j}$ denotes the partial derivative of $\xi$ with respect to $\hat{x}_{j}$ and $\hat{t}$ represents the time. Equations \eqref{le1}-\eqref{le3} represent mass, linear and angular momentum conservation, respectively, with $\Omega$ being the domain occupied by the liquid crystal, $\rho$ is the mass density (assumed constant) and $\rho_{1}$ is a constant, measured in terms of weight divided by distance, that arises from the consideration of a rotational kinetic energy of the material element. Here, $\sigma$, $\pi$ and $g$ represent, respectively, the stress tensor, the director stress tensor and the intrinsic director body force. They are defined as
\begin{equation}
\begin{array}{l l}
\sigma_{ji}= & -P\delta_{ij} - \displaystyle \frac{{\rm d}F}{{\rm d} n_{k,j}}n_{k,i}+\bar{\sigma}_{ji},\\
\pi_{ji}= & \beta_{j} n_{i} + \displaystyle \frac{{\rm d}F}{{\rm d}n_{i,j}},\\
g_{i}= &\gamma n_{i}-\beta_{j}n_{i,j} - \displaystyle \frac{ {\rm d} F}{{\rm d} n_{i}} +  \bar{g}_{i},
\end{array}
\end{equation}
where $P$ is the pressure of the fluid flow and $\delta_{ij}$ is the Kronecker delta. The vector $(\beta_{1},\beta_{2},\beta_{3})$ and the scalar function $\gamma$ (sometimes called direction tension) are Lagrange multipliers ensuring $\|\ve{n}\|=1$.\cite{LeslieGeneral} $F$ represents the Frank--Oseen free elastic energy, which is associated to distortions of the anisotropic axis. In the case of nematic liquid crystals, $F$ depends on four elastic constants $K_{i}$ ($i=1,2,3,4$), corresponding to the curvature components describing splay, twist, bend and saddle-splay effects (see for instance equation (4.130) in Stewart\cite{stewart2004static}). Here, we exploit the one-constant approximation of the Frank--Oseen elastic free energy density given by\cite{p1995physics}
$$ F =  \frac{K}{2} n_{i,j}n_{i,j},$$
where $K$ is the only elastic constant considered. Note that in this framework, $F$ does not depend on $n_{i}$, so that the term $\frac{{\rm d}F}{{\rm d}n_{i}}$ appearing in the definition of $g_{i}$ can be disregarded. Furthermore,
$$\bar{\sigma}_{ji}= \alpha_{1}n_{k}n_{p}A_{kp}n_{i}n_{j} + \alpha_{2}N_{i}n_{j} + \alpha_{3} N_{j}n_{i} + \alpha_{4}A_{ij} + \alpha_{5}A_{ik}n_{k}n_{j} + \alpha_{6}A_{jk}n_{k}n_{i},$$
$$\begin{array}{c} \displaystyle N_{i}= \frac{{\rm d} n_{i}}{{\rm d} \hat{t}}+ \ve{v}\cdot \nabla n_{i} - \omega_{i j}n_{j}, \mbox{ } \omega_{ij}= \displaystyle \frac{v_{i,j} - v_{j,i}}{2}, \mbox{ } A_{ij}= \displaystyle \frac{v_{i,j} + v_{j,i}}{2},\\
\\
\text{ and } \bar{g}_{i} = - \gamma_{1}N_{i} - \gamma_{2}n_{j}A_{ji},\end{array}$$
where $\alpha_{i}$ are constant viscosities satisfying the Parodi relation,\cite{Parodi} $\alpha_{2}+\alpha_{3}=\alpha_{6}-\alpha_{5}$, and $\gamma_{1}=\alpha_{3}-\alpha_{2}$, $\gamma_{2}=\alpha_{6}-\alpha_{5}$. More details about these parameters can be found in Section \ref{rmk:coef}.

\subsection{Simplified model}\label{sec:appendixsimplified}
In this work, we assume that the microfluidic channel, with domain $(0,l)\times(0,w)\times(-h,h)$ has length $l$ much greater than width $w$ and width much greater than height $2h$, so that the director and the flow fields may be assumed to depend only on the $\hat{z}$-coordinate. Thus, we let $\ve{n}=(\sin(\theta(\hat{z},\hat{t}),0,\cos(\theta(\hat{z},\hat{t}))$, $\ve{v}=(u(\hat{z},\hat{t}),0,0)$. Moreover, $u(\hat{z},\hat{t})$ is considered symmetric around $\hat{z}=0$ and the no-slip condition is assumed in the channel walls (i.e. $u(\pm h,\hat{t})=0$).
\\
Using this information in the constitutive formulae, one has that\\
\\
$\bullet$ $A_{ij}=0$ except for $A_{13}=A_{31}= \displaystyle\frac{u_{\hat{z}}}{2}$.\\
$\bullet$ $\omega_{ij}=0$ except for $\omega_{13}=\displaystyle\frac{u_{\hat{z}}}{2}$ and $\omega_{31}= \displaystyle\frac{-u_{\hat{z}}}{2}$.\\
$\bullet$ $N_{1}=n_{1,\hat{t}}- w_{13}n_{3}= \cos(\theta)\theta_{\hat{t}} - \displaystyle\frac{u_{\hat{z}}}{2} \cos(\theta) = \cos(\theta)(\theta_{\hat{t}} - \displaystyle\frac{u_{\hat{z}}}{2})$.\\ 
$\bullet$ $N_{2}=0$.\\
$\bullet$ $N_{3}= n_{3,\hat{t}} - w_{31}n_{1}= -\sin(\theta)\theta_{\hat{t}} + \displaystyle \frac{u_{\hat{z}}}{2} \sin(\theta) = \sin(\theta) (\displaystyle\frac{u_{\hat{z}}}{2}-\theta_{\hat{t}})$.\\
$\bullet$ $\bar{g}_{1}= - \gamma_{1} N_{1} - \gamma_{2}A_{31}n_{3}= \displaystyle \frac{\cos(\theta)u_{\hat{z}}}{2}(\gamma_{1}-\gamma_{2}) - \gamma_{1}\cos(\theta)\theta_{\hat{t}}$.\\
$\bullet$ $\bar{g}_{2}=0$.\\
$\bullet$ $\bar{g}_{3}= - \gamma_{1} N_{3} - \gamma_{2}A_{13}n_{1}= -\displaystyle \frac{\sin(\theta)u_{\hat{z}}}{2}(\gamma_{1}+\gamma_{2}) +\gamma_{1}\sin(\theta)\theta_{\hat{t}}.$\\

Now, taking into account that $F$ only depends on the variables $n_{1,3}$ and $n_{3,3}$ one has that $\pi_{ij,i}=0$ except for $\pi_{31,3}$ and $\pi_{33,3}$. Thus,\\
\\
$\bullet$ $\pi_{31,3}= \displaystyle \left (\frac{{\rm d}F}{{\rm d} n_{1,3}}\right)_{,3}= Kn_{1,33}.$\\
$\bullet$ $\pi_{33,3}= \displaystyle \left (\frac{{\rm d}F}{{\rm d} n_{3,3}}\right)_{,3}= Kn_{3,33}.$\\
$\bullet$ $g_{1}= \gamma n_{1} + \bar{g}_{1}= \gamma n_{1} +  \displaystyle \frac{\cos(\theta)u_{\hat{z}}}{2}(\gamma_{1}-\gamma_{2}) - \gamma_{1}\cos(\theta)\theta_{\hat{t}}.$\\
$\bullet$ $g_{2}= 0.$\\
\\
$\bullet$ $g_{3}= \gamma n_{3} + \bar{g}_{3}= \gamma n_{3} - \displaystyle \frac{\sin(\theta)u_{\hat{z}}}{2}(\gamma_{1}+\gamma_{2}) +\gamma_{1}\sin(\theta)\theta_{\hat{t}}.$\\
$\bullet$ $\bar{\sigma}_{ij}=0 \mbox{ except for } \bar{\sigma}_{11},\bar{\sigma}_{13},\bar{\sigma}_{31} \mbox{ and } \bar{\sigma}_{33}.$\\

In our case, it follows from the linear momentum equation (\ref{le2}) that 
\begin{equation*}
\begin{array}{r l}
\rho\frac{{\rm d}u}{{\rm d} \hat{t}} & = \sigma_{11,1}+\sigma_{31,3}= -P_{,1} + \bar{\sigma}_{31,3},
\\
\\
0  & = \sigma_{22,2} = - P_{,2} ,
\\
\\
0 & = \sigma_{33,3} = -(P+ 2F)_{,3} + \bar{\sigma}_{33,3}.
\end{array}
\end{equation*}
Note that we will use the notation $f_{,1}$, $f_{,2}$, $f_{,3}$; $f_{,\hat{x}}$, $f_{,\hat{y}}$, $f_{,\hat{z}}$ and $f_{\hat{x}}$, $f_{\hat{y}}$, $f_{\hat{z}}$ interchangeably. Therefore, it follows from (\ref{le2}) that
\begin{subequations}
\begin{alignat}{2}
-(P+2F)_{\hat{x}} + \bar{\sigma}_{31,\hat{z}} & = \rho\frac{{\rm d}u}{{\rm d} \hat{t}}&&\mbox{ in }(0,l)\times(0,w)\times(-h,h)\times(0,+\infty),\label{i1}
\\
 (P+ 2F)_{\hat{y}} &= 0&&\mbox{ in }(0,l)\times(0,w)\times(-h,h)\times(0,+\infty),\label{i2}
\\
-(P+ 2F)_{\hat{z}} + \bar{\sigma}_{33,\hat{z}} & = 0 && \mbox{ in }(0,l)\times(0,w)\times(-h,h)\times(0,+\infty).\label{i3}
\end{alignat}
\end{subequations}
We suppose that the inertia of the liquid crystal molecules can be ignored in typical cells having small depths,\cite{stewart2004static} so that $\rho \displaystyle \frac{{\rm d}u}{{\rm d} \hat{t}} =0$ in equation (\ref{i1}).
From (\ref{i2}), one has that $P+ 2F= q(\hat{x},\hat{z},\hat{t})$. Now, if we integrate with respect to $\hat{x}$ in equation (\ref{i1}) and take into account that $F$ only depends on $\hat{z}$ and $\hat{t}$,
\begin{equation}\label{i4}
P+ 2F= \hat{x} \bar{\sigma}_{31,\hat{z}} + r(\hat{z},\hat{t}).
\end{equation}
If relation (\ref{i4}) is introduced in equation (\ref{i3}), one has that $(\hat{x} \bar{\sigma}_{31,\hat{z}} + r(\hat{z},\hat{t}))_{,\hat{z}}= \bar{\sigma}_{33,\hat{z}}$. Consequently, $\bar{\sigma}_{31,\hat{z}\hat{z}}=0,$ and so
\begin{equation}\label{eqsigma}
\bar{\sigma}_{31}=C(\hat{t})\hat{z}+D(\hat{t}),
\end{equation}
where $C(\hat{t})$ and $D(\hat{t})$ are functions to be determined. Then, from relation (\ref{i4}), one has that
\begin{equation}\label{i5}
P+ 2F= C(\hat{t})\hat{x} + r(\hat{z},\hat{t}).
\end{equation}
From equations (\ref{i3}) and (\ref{i5}) it follows that $(C(\hat{t})\hat{x}+r(\hat{z},\hat{t}))_{,\hat{z}}= (r(\hat{z},\hat{t}))_{,\hat{z}}=  \bar{\sigma}_{33,\hat{z}}$, where integrating with respect to $\hat{z}$ one has that $r(\hat{z},\hat{t})=  \bar{\sigma}_{33} + s(\hat{t})$, $s$ being a time-dependent function to be determined. Returning to equation (\ref{i5}), it follows that
\begin{equation}\label{eqp}
P= - 2F + C(\hat{t})\hat{x} + s(\hat{t}) + \bar{\sigma}_{33}.
\end{equation}
Replacing the value of $\bar{\sigma}_{31}$ in equation (\ref{eqsigma}) one has that
$$u_{\hat{z}} g(\theta) + \theta_{\hat{t}} m(\theta) = C(\hat{t})\hat{z}+D(\hat{t}).$$

A consequence of the symmetry of $u$ enforces $\frac{\partial \theta}{\partial \hat{t}} = 0$ at $\hat{z}=0$. Any scenario for which $\frac{\partial \theta}{\partial \hat{t} } \neq 0$ would induce a velocity profile that is non-symmetric and thus violate our original assumption. As a result, this implies that $D(\hat{t})=0$ for our system and hence
\begin{equation}\label{deduccion1}
u_{\hat{z}} g(\theta) + \theta_{\hat{t}} m(\theta) = C(\hat{t})\hat{z},
\end{equation}
where
\begin{subequations}
\begin{align}
g(\theta) & = \alpha_{1}\cos^{2}(\theta)\sin^{2}(\theta)  +\frac{\alpha_{5}-\alpha_{2}}{2}\cos^{2}(\theta) +  \frac{\alpha_{3}+\alpha_{6}}{2} \sin^{2}(\theta)+ \frac{\alpha_{4}}{2},\\
m(\theta)& =  \alpha_{2}\cos^{2}(\theta)-\alpha_{3}\sin^{2}(\theta).
\end{align}
\end{subequations}
Note that we have reduced equations (\ref{i1})--(\ref{i3}) to equation (\ref{deduccion1}), the pressure being available via equation (\ref{eqp}).
Now, the angular momentum balance equation (\ref{le3}) for $i=1$ and $i=3$ reduces, respectively, to
$$\begin{array}{l l } \rho_{1} n_{1,\hat{t}\hat{t}} & = g_{1} + \pi_{31,3} = \gamma n_{1} + \bar{g}_{1} + \pi_{31,3}= \gamma n_{1} + \bar{g}_{1} + K n_{1, 33}, \\
\\
 \rho_{1} n_{3,\hat{t}\hat{t}} & = g_{3} + \pi_{33,3} = \gamma n_{3} + \bar{g}_{3} + \pi_{33,3}= \gamma n_{3} + \bar{g}_{3} + K n_{3, 33}.
\end{array}
$$
It remains to compute $n_{1,33}$, $n_{3,33}$, $n_{1,\hat{t}\hat{t}}$ and $n_{2,\hat{t}\hat{t}}$:
\begin{itemize}
\item $n_{1}= \sin(\theta) \Rightarrow n_{1,3}=\cos(\theta)\theta_{\hat{z}} \Rightarrow n_{1,33}= -\sin(\theta)(\theta_{\hat{z}})^2 + \cos(\theta)\theta_{\hat{z}\hat{z}}$,
\item  $n_{1,\hat{t}}=\cos(\theta)\theta_{\hat{t}} \Rightarrow n_{1,\hat{t}\hat{t}}= -\sin(\theta)(\theta_{\hat{t}})^2 + \cos(\theta)\theta_{\hat{t}\hat{t}}$,
\item $n_{3}=\cos(\theta)  \Rightarrow n_{3,3}=-\sin(\theta)\theta_{\hat{z}} \Rightarrow n_{3,33}= -\cos(\theta)(\theta_{\hat{z}})^2 - \sin(\theta)\theta_{\hat{z}\hat{z}}$,
\item $n_{3,\hat{t}}=-\sin(\theta)\theta_{\hat{t}} \Rightarrow n_{3,\hat{t}\hat{t}}= -\cos(\theta)(\theta_{\hat{t}})^2 - \sin(\theta)\theta_{\hat{t}\hat{t}}.$
\end{itemize}
Thus, equation (\ref{le3}) when $i=1$ and $i=3$ becomes
$$\begin{array}{r l} \rho_{1} (- \sin(\theta)(\theta_{\hat{t}})^2 + \cos(\theta)\theta_{\hat{t}\hat{t}}) = & \gamma \sin(\theta) - \gamma_{1} \cos(\theta) \theta_{\hat{t}} + \cos(\theta)\displaystyle \frac{u_{\hat{z}}}{2}(\gamma_{1}-\gamma_{2})  \\
& + K(-\sin(\theta)\theta_{\hat{z}}^{2} + \cos(\theta)\theta_{\hat{z}\hat{z}}),\\
\\
\rho_{1} (- \cos(\theta)(\theta_{\hat{t}})^2 - \sin(\theta)\theta_{\hat{t}\hat{t}}) = & \gamma \cos(\theta) + \gamma_{1} \sin(\theta) \theta_{\hat{t}} - \sin(\theta)\displaystyle \frac{u_{\hat{z}}}{2}(\gamma_{1}+\gamma_{2})  \\
&+ K(-\cos(\theta)\theta_{\hat{z}}^{2} - \sin(\theta)\theta_{\hat{z}\hat{z}}).
\end{array}$$
We neglect the term $\rho_{1}\theta_{\hat{t}\hat{t}}$ (it is accepted as being negligible in physical situations\cite{stewart2004static}). Then, multiplying the first equation by $\cos(\theta)$, the second one by $\sin(\theta)$ and subtracting them, one obtains:
\begin{equation}\label{deduccion2}
\gamma_{1}\theta_{\hat{t}} = K \theta_{\hat{z}\hat{z}} + \frac{u_{\hat{z}}}{2} \big(\gamma_{1}-\gamma_{2}\cos(2\theta)\big).
\end{equation}
Thus, the evolution of $\theta$ and $u$ are described by the following system
\begin{subequations}\label{ledynamic}
\begin{align}
\gamma_{1}\theta_{\hat{t}} =  & K \theta_{\hat{z}\hat{z}} - u_{\hat{z}} m(\theta) & \hat{z}\in (-h,h), \hat{t}>0, \label{ledynamic1}
\\
 C(\hat{t})\hat{z}  =  & u_{\hat{z}} g(\theta) + \theta_{\hat{t}} m(\theta) & \hat{z}\in (-h,h), \hat{t}>0, \label{ledynamic2}
\\
\theta(\hat{z},0) = & \Theta(\hat{z}) & \hat{z} \in (-h,h),\label{ledynamic3}
\\
u(\pm h,\hat{t}) =& 0 &  \hat{t}>0,\label{ledynamic4}
\end{align}
\end{subequations}
where $\Theta$ is the initial condition for $\theta$ and $C(\hat{t})=P_{\hat{x}}$, i.e, the channel direction component of the pressure gradient.
\subsubsection{Remarks on coefficients}\label{rmk:coef}
The coefficients $\alpha_{i}$ and $\gamma_{i}$ are usually called \textit{Leslie Coefficients} (see for instance Lee\cite{Lee1988} or Wang et al.\cite{Wang2006} for further information about their physical meaning and how to approximate them experimentally). They depend only on the temperature and have the dimension of viscosity. Some constraints on the \textit{Leslie Coefficients} come from the non-negativity of the \textit{Dissipative function}.\cite{Leslie1968,stewart2004static} When the Parodi relation is used\cite{Parodi}, the dissipative function is defined as:\cite{stewart2004static}
$$\mathcal{D}= \alpha_{1}(n_{i}A_{ij}n_{j})^2 + 2\gamma_{2}N_{i}A_{ij}n_{j} + \alpha_{4}A_{ij}A_{ij} + (\alpha_{5}+\alpha_{6})A_{ij}A_{jk}n_{i}n_{k} + \gamma_{1}N_{i}N_{i}.$$
In our particular case,
$$\begin{array}{r l}\mathcal{D}  = & \alpha_{1} u_{\hat{z}}^{2}\sin^{2}(\theta)\cos^{2}(\theta)+ 2\gamma_{2}\displaystyle \frac{u_{\hat{z}}}{2}(\theta_{\hat{t}}-\frac{u_{\hat{z}}}{2})(\cos^{2}(\theta)-\sin^{2}(\theta)) + \alpha_{4}\displaystyle\frac{u_{\hat{z}}^{2}}{2} \\
\\
&+ (\alpha_{5}+\alpha_{6})\displaystyle\frac{u_{\hat{z}}^{2}}{4} + \gamma_{1}(\theta_{\hat{t}}-\displaystyle\frac{u_{\hat{z}}}{2})^2 = 2\theta_{\hat{t}}u_{\hat{z}}m(\theta) + \gamma_{1}\theta_{\hat{t}}^{2} + g(\theta)u_{\hat{z}}^{2}.\end{array}$$
This expression is a quadratic form and can be rewritten as:
$$\mathcal{D}= \left [\begin{array}{c c} X & Y \end{array} \right] \left [\begin{array}{c c} g(\theta) & m(\theta) \\ m(\theta) & \gamma_{1} \end{array}\right]\left[\begin{array}{c} X \\ Y \end{array} \right], \mbox{ with } X= u_{\hat{z}} \mbox{, } Y=\theta_{\hat{t}}.$$
A reasonable assumption is that the dissipation function is positive,\cite{stewart2004static} which is fulfilled if and only if the determinant of every principal submatrix is positive,\cite{Infante2015} i.e., 
\begin{equation}\label{inequalities}
g(\theta)> 0 \hspace{1cm} \mbox{ and }\hspace{1cm} \gamma_{1}g(\theta)-m^2(\theta)> 0.
\end{equation}
When $\theta=0$, this implies that
$$\gamma_{1}> 0, \hspace{5mm} \alpha_{5}-\alpha_{2}+\alpha_{4}> 0 \hspace{5mm}\mbox{and}\hspace{5mm}  \gamma_{1}(\alpha_{5}-\alpha_{2}+\alpha_{4})> 2\alpha_{2}^{2}.$$ 
\section{Equilibrium Solutions with $\mathcal{G}=0$.}\label{apsec:steadyg0}
\begin{figure}[ht!]
\centering
%
%
\begin{psfrags}%
\psfragscanon%
%
\psfrag{s02}[lt][lt]{\color[rgb]{0,0,0}\setlength{\tabcolsep}{0pt}\begin{tabular}{l}{\Large a}\end{tabular}}%
\psfrag{s03}[rt][rt]{\color[rgb]{0,0,0}\setlength{\tabcolsep}{0pt}\begin{tabular}{r}{\Large b}\end{tabular}}%
\psfrag{s04}[b][b]{\color[rgb]{0,0,0}\setlength{\tabcolsep}{0pt}\begin{tabular}{c}{\Large $\mathcal{B}$}\end{tabular}}%
%
\psfrag{x01}[t][t]{-$2\pi$}%
\psfrag{x02}[t][t]{-$\pi$}%
\psfrag{x03}[t][t]{0}%
\psfrag{x04}[t][t]{$\pi$}%
\psfrag{x05}[t][t]{$2\pi$}%
%
\psfrag{v01}[r][r]{0}%
\psfrag{v02}[r][r]{$\frac{\pi}{4}$}%
\psfrag{v03}[r][r]{$\frac{\pi}{2}$}%
%
\psfrag{z01}[r][r]{0}%
\psfrag{z02}[r][r]{0.5}%
\psfrag{z03}[r][r]{1}%
\psfrag{z04}[r][r]{1.5}%
\psfrag{z05}[r][r]{2}%
\psfrag{z06}[r][r]{2.5}%
%
\resizebox{10cm}{!}{\includegraphics{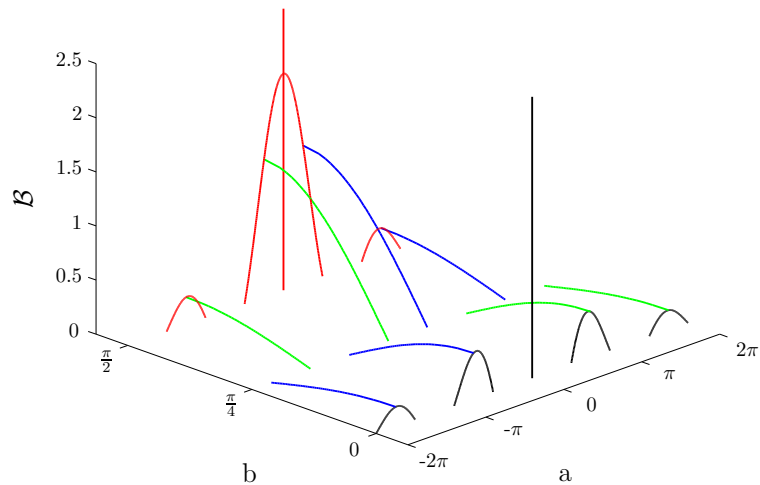}}%
\end{psfrags}%
%

\caption{Solution landscape with $a\in [-2\pi, 2\pi]$ and $b\in[0,\frac{\pi}{2}]$. Solutions of Type I (Type II) correspond to $b=0$ ($b=\frac{\pi}{2}$) and are plotted in black (red). Solutions of Types III and IV correspond to $b\in [0,\frac{\pi}{2}]$ and are plotted in blue and green, respectively.}
\label{fig:landscapeg03d}
\end{figure}

When $\mathcal{G}=0$, we can explicitly solve the first equation
of system (\ref{systemstatic}) to obtain $\theta^{\ast}(z)= a z + b$ where $a$ and $b$ are constants to be determined by the boundary conditions. These solutions may be categorized as  
\begin{align}
\mbox{ Type I } & \theta^{\ast}(z) = a_{n}z + m\pi, & \mbox{ where }  m \in \mathds{Z} \mbox{ and }\mathcal{B}a_{n}= -\sin(2a_{n}),  \label{eq:type1}
\\ 
\mbox{ Type II } & \theta^{\ast}(z) = \tilde{a}_{n}z + (m + \frac{1}{2})\pi, & \mbox{ where }m \in \mathds{Z} \mbox{ and }  \mathcal{B}\tilde{a}_{n}= \sin(2\tilde{a}_{n}), \label{eq:type2} 
\\
\mbox{ Type III } & \theta^{\ast}(z) = (n+\frac{1}{4})\pi z + b_{m}, & \mbox{ where } n \in \mathds{Z} \mbox{ and } \nonumber
\\
& & \cos(2 b_{m})= -\mathcal{B}(n+\frac{1}{4})\pi, \label{eq:type3} 
\\
\mbox{ Type IV } & \theta^{\ast}(z) = (n+\frac{3}{4})\pi z + \tilde{b}_{m}, & \mbox{ where } n \in \mathds{Z} \mbox{ and } \nonumber \\
& & \cos(2 \tilde{b}_{m})= \mathcal{B}(n+\frac{3}{4})\pi. \label{eq:type4}
\end{align}
For every value of $\mathcal{B}$, we obtain ordered set of solutions for \eqref{eq:type1}, with $0=a_{0}<a_{1}<\ldots<a_{n}$ ($n \in \mathds{N}\cup \{0\}$ depending on $\mathcal{B}$). Furthermore, if $a_{n}$ provides a solution, so does $-a_{n}$, which we denote by $a_{-n}$. Equivalent statement can be made for $\tilde{a}_{n}$, $b_{m}$ and $\tilde{b}_{m}$, solutions of equations \eqref{eq:type2}, \eqref{eq:type3} and \eqref{eq:type4}, respectively. \\
We observe that constant solutions of Type I and II, $\theta^{\ast} \equiv k\frac{\pi}{2}$ ($k\in \mathds{Z}$) exist for all values of $\mathcal{B}$, while solutions of Type III and IV exist only if $\mathcal{B}\leq\frac{4}{\pi}$. The associated director fields are
$$\begin{array}{r l}
\mbox{ Type I 	 }& \ve{n}(z)= (-1)^{m} (\sin(a_{n}z),0,\cos(a_{n}z)),
\\
\mbox{ Type II  	} & \ve{n}(z) = (-1)^{m}(\cos(\tilde{a}_{n}z),0,-\sin(\tilde{a}_{n}z)),
\\
\mbox{ Type III  	} & \ve{n}(z) = (-1)^{m}(\sin((n+\frac{1}{4})\pi z+ b_{0}),0,\cos((n+\frac{1}{4})\pi z+ b_{0}))
\\
\mbox{ Type IV  	} & \ve{n}(z)= (-1)^{m}(\sin((n+\frac{3}{4})\pi z+ \tilde{b}_{0}),0,\cos((n+\frac{3}{4})\pi z+ \tilde{b}_{0}))
\end{array}$$
and thus, since director fields with $m\in \mathds{Z}$ are the same but with opposite direction, all possible director profiles in \eqref{eq:type1}--\eqref{eq:type4} are covered by $m=0$.
Figure \ref{fig:landscapeg03d} shows the solution landscape in terms of $a$, $b$ and $\mathcal{B}$, restricted to $a\in [2\pi,2\pi]$ and $b \in [0,\frac{\pi}{2}]$.
Since solutions of Types III and IV are always unstable (see Section \ref{apsec:linearstability}), we only track solutions of Type I and II in this paper. 
\subsection{Linear Stability of Equilibrium Solutions}\label{apsec:linearstability}
We analyze the linear stability of the equilibria \eqref{eq:type1}--\eqref{eq:type4} by linearizing around the steady state, so that $\theta(z,t)
\approx \theta^{\ast}(z)+\bar{\theta}(z,t)$, with $\bar{\theta}(z,t)$ satisfying
\begin{equation}\label{linearized}
\left \{\begin{array}{l} \bar{\theta}_{t}=F(\theta^{\ast}(z))\bar{\theta}_{zz}, \\ \bar{\theta}(z,0)= \delta \Theta(z),\\
\mathcal{B}\bar{\theta}_{z}(1,t)=-2\cos(2\theta^{\ast}(1))\bar{\theta}(1,t),\\
\mathcal{B}\bar{\theta}_{z}(-1,t)=2\cos(2\theta^{\ast}(-1))\bar{\theta}(-1,t), \end{array}\right.
\end{equation}
where $F(\theta^{\ast})= \displaystyle \frac{g(\theta^{\ast})}{\gamma_{1}g(\theta^{\ast}) + m(\theta^{\ast})h(\theta^{\ast})}$ and $\delta \Theta(z)$ being a small perturbation of $\theta^{\ast}$. It is straightforward to show that (\ref{linearized}) admits a separable solution of the form 
$$\bar{\theta}(z,t) =\sum_{k=0}^{\infty} C_{k}e^{-\lambda_{k}t}Z_{k}(z),$$
for suitable eigenvalues $\{\lambda_{k}\}_{k\in \mathds{N}}\subset \mathds{R}$ and $\{C_{k}\}_{k\in \mathds{N}} \subset \mathds{R}$ such that $\delta \Theta(z) = \sum_{k=0}^{\infty} C_{k}Z_{k}(z)$, $\{Z_{k}\}_{k\in \mathds{N}}$ solving the following second-order ordinary differential equation
\begin{equation}
\begin{array}{l} F(\theta^{\ast}(z))Z_{k}''(z) + \lambda_{k}Z_{k}(z)=0 \\ \mathcal{B}Z_{k}'(1)=-2\cos(2\theta^{\ast}(1))Z_{k}(1)  \\ \mathcal{B}Z_{k}'(-1)=2\cos(2\theta^{\ast}(-1))Z_{k}(-1). \end{array} 
\end{equation}
When $\theta^{\ast}(z)=0$, we find that $\lambda_{k}$ must be positive in order to find a solution of \eqref{linearized}. Particularly, 
\begin{align}
Z_{k}(z)=\sin\left(\sqrt{\frac{\lambda_{k}}{F(0)}}\frac{z+1}{2}\right)+ \frac{\mathcal{B}}{2}\sqrt{\frac{\lambda_{k}}{F(0)}}
\cos\left(\sqrt{\frac{\lambda_{k}}{F(0)}}\frac{z+1}{2}\right),
\end{align}
where
\begin{align}
F(0)= \frac{\alpha_{5}-\alpha_{2}+1}{\gamma_{1}(\alpha_{5}-\alpha_{2}+1)-2\alpha_{2}^{2}} >0
\end{align}
using (\ref{strictinequalities}), and $\lambda_{k}$ satisfies the transcendental equation
\begin{equation}\label{lambdatheta0}
\tan\left( \sqrt{\frac{\lambda_{k}}{F(0)}}\right) =
-\frac{4\mathcal{B}\sqrt{\frac{\lambda_{k}}{F(0)}}}{4-\mathcal{B}^{2}\frac{\lambda_{k}}{F(0)}}.
\end{equation}
Thus, since $\lambda_{k}>0$ $\forall$ $k=0, 1,\ldots$, $\bar{\theta}(z,t) \xrightarrow{t\rightarrow \infty} 0$ and  the trivial solution $\theta^{\ast}=0$ is linearly stable in
this dynamic framework. For the other steady states in \eqref{eq:type1}--\eqref{eq:type4}, we follow the same paradigm as above and
numerically compute the eigenvalues $\lambda_{k}$ using the function \textit{eigs}, in the MATLAB package \textit{Chebfun}
(\url{http://www.chebfun.org}). 
We find that in terms of the type of solution $\theta^{\ast}$, the stability can be classified as:
$$\begin{array}{c c}
\mbox{Type I 		} & \mbox{ is stable if }n \mbox{ is even and unstable if }n \mbox{ is odd }\\
\mbox{Type II 		} & \mbox{ is stable if }n \mbox{ is odd and unstable if }n \mbox{ is even }\\
\mbox{Type III - IV 		} &  \mbox{ is unstable. }
\end{array}$$
\subsubsection{Sample liquid crystal molecular configurations}\label{sec:config}
\begin{figure}[ht!]
\centering
\subfigure[$\theta^{\ast}_{a_{-4}}$]{\includegraphics[scale=0.33]{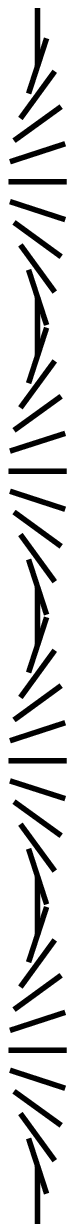}}
\subfigure[$\theta^{\ast}_{a_{-3}}$]{\includegraphics[scale=0.33]{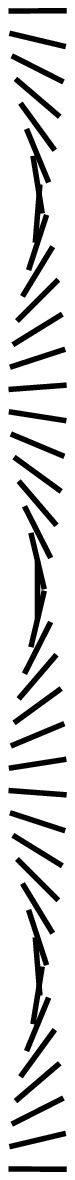}}
\subfigure[$\theta^{\ast}_{a_{-2}}$]{\includegraphics[scale=0.33]{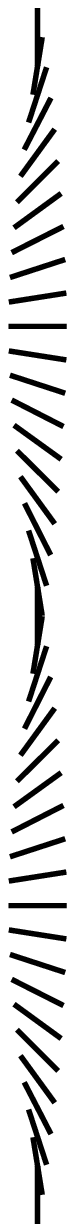}}
\subfigure[$\theta^{\ast}_{a_{-1}}$]{\includegraphics[scale=0.33]{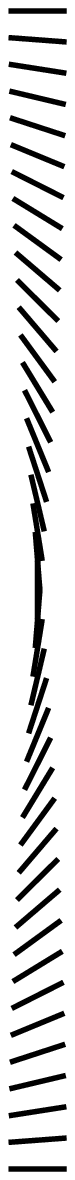}}
\subfigure[$\theta^{\ast}_{a_{0}}$]{\includegraphics[scale=0.33]{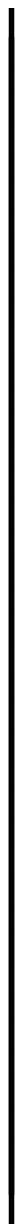}}
\subfigure[$\theta^{\ast}_{a_{1}}$]{\includegraphics[scale=0.33]{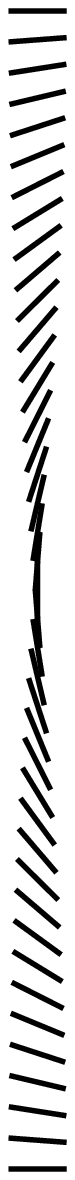}}
\subfigure[$\theta^{\ast}_{a_{2}}$]{\includegraphics[scale=0.33]{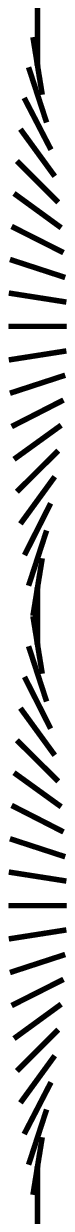}}
\subfigure[$\theta^{\ast}_{a_{3}}$]{\includegraphics[scale=0.33]{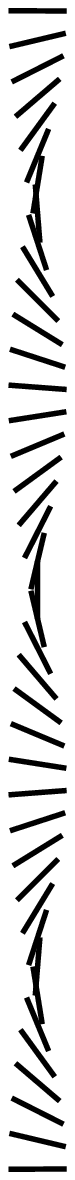}}
\subfigure[$\theta^{\ast}_{a_{4}}$]{\includegraphics[scale=0.33]{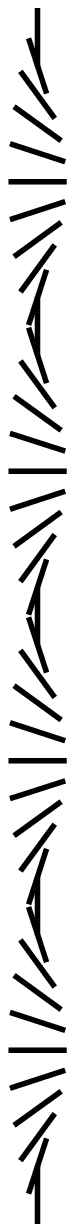}}
\caption{$\ve{n}$ associated with steady states $\theta^{\ast}_{a_{n}}$ (Type I), obtained with $\mathcal{B}=0.001$ and $\mathcal{G}=0$. These states are stable if $n$ is even and unstable if $n$ is odd.}
\label{fig:configtype1}
\end{figure}

\begin{figure}[ht!]
\centering
\subfigure[$\theta^{\ast}_{\tilde{a}_{-4}}$]{\includegraphics[scale=0.33]{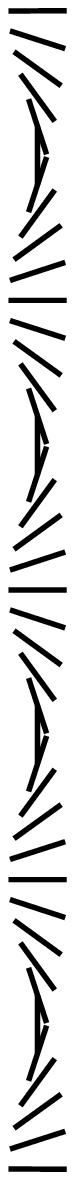}}
\subfigure[$\theta^{\ast}_{\tilde{a}_{-3}}$]{\includegraphics[scale=0.33]{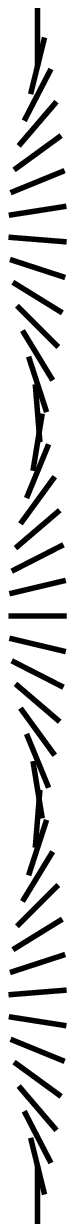}}
\subfigure[$\theta^{\ast}_{\tilde{a}_{-2}}$]{\includegraphics[scale=0.33]{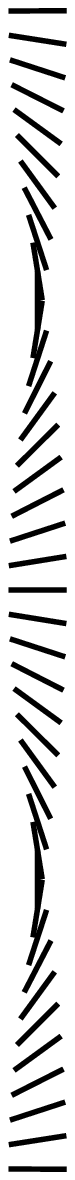}}
\subfigure[$\theta^{\ast}_{\tilde{a}_{-1}}$]{\includegraphics[scale=0.33]{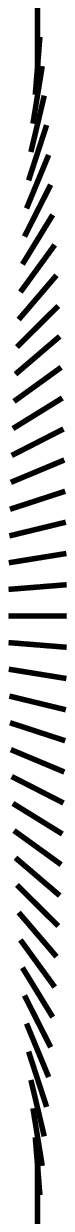}}
\subfigure[$\theta^{\ast}_{\tilde{a}_{0}}$]{\includegraphics[scale=0.33]{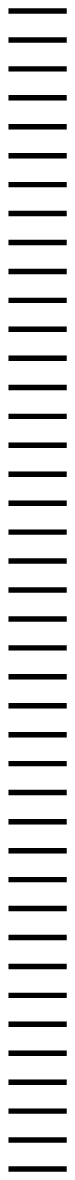}}
\subfigure[$\theta^{\ast}_{\tilde{a}_{1}}$]{\includegraphics[scale=0.33]{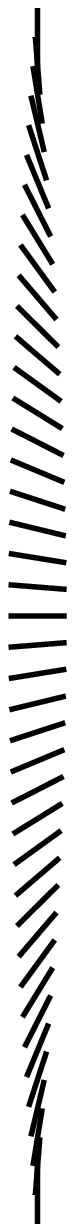}}
\subfigure[$\theta^{\ast}_{\tilde{a}_{2}}$]{\includegraphics[scale=0.33]{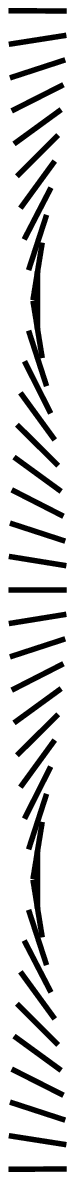}}
\subfigure[$\theta^{\ast}_{\tilde{a}_{3},0}$]{\includegraphics[scale=0.33]{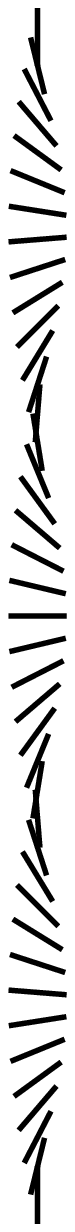}}
\subfigure[$\theta^{\ast}_{\tilde{a}_{4},0}$]{\includegraphics[scale=0.33]{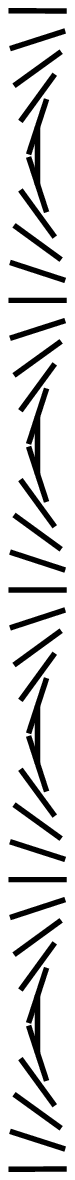}}
\caption{$\ve{n}$ associated with steady states $\theta^{\ast}_{\tilde{a}_{n}}$ (Type II), obtained with $\mathcal{B}=0.001$ and $\mathcal{G}=0$. These states are stable if $n$ is odd and unstable if $n$ is even.}
\label{fig:configtype2}
\end{figure}

In this section we show the director field corresponding to some steady state solutions of the system (\ref{systemtdep}) with $\mathcal{G}=0$. Particularly, Figures \ref{fig:configtype1} and \ref{fig:configtype2} show the director profiles associated, respectively, to solutions $\theta^{\ast}_{a_{n}}$ and $\theta^{\ast}_{\tilde{a}_{n}}$, $n=0, \pm 1, \pm2 \pm 3, \pm 4$.
\newpage
\bibliographystyle{plain}
\bibliography{bibliography}

\begin{thebibliography}{10}

\bibitem{allgower2003}
E.~L. Allgower and K.~Georg.
\newblock {\em Introduction to Numerical Continuation Methods}.
\newblock Classics in Applied Mathematics. Society for Industrial and Applied
  Mathematics, 2003.

\bibitem{Linda2015}
T.~G. Anderson, E.~Mema, L.~Kondic, and L.~J. Cummings.
\newblock Transitions in {P}oiseuille flow of nematic liquid crystal.
\newblock {\em International Journal of Non-Linear Mechanics}, 75:15 -- 21,
  2015.

\bibitem{dg}
P.~G. de~Gennes and J.~Prost.
\newblock {\em The Physics of Liquid Crystals}.
\newblock International Series of Monographs on Physics. Clarendon Press, 1993.

\bibitem{p1995physics}
P.~G. de~Gennes and J.~Prost.
\newblock {\em The Physics of Liquid Crystals}.
\newblock International Series of Monographs on Physics. Clarendon Press, 1995.

\bibitem{yeomans}
C.~Denniston, E.~Orlandini, and J.~M. Yeomans.
\newblock Lattice {B}oltzmann simulations of liquid crystal hydrodynamics.
\newblock {\em Phys. Rev. E}, 63:056702, Apr 2001.

\bibitem{Infante2015}
J.~A. Infante and J.~M. Rey.
\newblock {\em M\'etodos Num\'ericos. Teor\'ia, problemas y pr\'acticas con
  MATLAB.}
\newblock Ediciones Pir\'amide. Grupo Anaya, 2015.

\bibitem{fancyref}
Y.~K. Kim, B.~Senyuk, and O.~D. Lavrentovich.
\newblock Molecular reorientation of a nematic liquid crystal by thermal
  expansion.
\newblock {\em Nature Communications}, 3(1133), 2012.

\bibitem{microfluidics}
O.~D. Lavrentovich.
\newblock Transport of particles in liquid crystals.
\newblock {\em Soft Matter}, 10:1264--1283, 2014.

\bibitem{Lee1988}
S.~Lee.
\newblock The {L}eslie coefficients for a polymer nematic liquid crystal.
\newblock {\em The Journal of Chemical Physics}, 88(8):5196--5201, 1988.

\bibitem{Leslie1968}
F.~M. Leslie.
\newblock Some constitutive equations for liquid crystals.
\newblock {\em Archive for Rational Mechanics and Analysis}, 28(4):265--283,
  1968.

\bibitem{ericksen}
F.~M. Leslie.
\newblock Continuum theory for nematic liquid crystals.
\newblock {\em Continuum Mechanics and Thermodynamics}, 4(3):167--175, 1992.

\bibitem{LeslieExistence}
F.~H. Lin and C.~Liu.
\newblock Existence of solutions for the {E}ricksen--{L}eslie system.
\newblock {\em Archive for Rational Mechanics and Analysis}, 154:135--156,
  2000.

\bibitem{mermin1979}
N.~D. Mermin.
\newblock The topological theory of defects in ordered media.
\newblock {\em Rev. Mod. Phys.}, 51:591--648, 1979.

\bibitem{Parodi}
O.~Parodi.
\newblock Stress tensor for a nematic liquid crystal.
\newblock {\em J. Phys. France}, 31(7):581--584, 1970.

\bibitem{Rapini1969}
A.~Rapini and M.~Papoular.
\newblock Distorsion d{'}une lamelle n\'ematique sous champ magn\'etique
  conditions d{'}ancrage aux parois.
\newblock {\em J. Phys. Colloques}, 30(C4), 1969.

\bibitem{senguptathesis}
A.~Sengupta.
\newblock {\em Topological Microfluidics: Nematic Liquid Crystals and Nematic
  Colloids in Microfluidic Environment}.
\newblock Springer Theses. Springer International Publishing, 2013.

\bibitem{Sengupta2013}
A.~Sengupta, U.~Tkalec, M.~Ravnik, J.~M. Yeomans, C.~Bahr, and S.~Herminghaus.
\newblock Liquid crystal microfluidics for tunable flow shaping.
\newblock {\em Phys. Rev. Lett.}, 110, Jan 2013.

\bibitem{sluckin2004crystals}
T.~J. Sluckin, D.~A. Dunmur, and H.~Stegemeyer.
\newblock {\em Crystals That Flow: Classic Papers from the History of Liquid
  Crystals}.
\newblock Liquid Crystals Book Series. Taylor \& Francis, 2004.

\bibitem{stewart2004static}
I.~W. Stewart.
\newblock {\em The Static and Dynamic Continuum Theory of Liquid Crystals: A
  Mathematical Introduction}.
\newblock Liquid Crystals Book Series. Taylor \& Francis, 2004.

\bibitem{stone2004engineering}
H.~A. Stone, A.~D. Stroock, and A.~Ajdari.
\newblock Engineering flows in small devices: microfluidics toward a
  lab-on-a-chip.
\newblock {\em Annu. Rev. Fluid Mech.}, 36:381--411, 2004.

\bibitem{Wang2006}
H.~Wang, T.~X. Wu, S.~Gauza, J.~R. Wu, and S.~Wu.
\newblock A method to estimate the {L}eslie coefficients of liquid crystals
  based on mbba data.
\newblock {\em Liquid Crystals}, 33(1):91--98, 2006.

\bibitem{whitesides2006}
G.~M. Whitesides.
\newblock The origins and the future of microfluidics.
\newblock {\em Nature}, 442(7101):368--373, 2006.

\bibitem{LeslieGeneral}
H.~Wu, X.~Xu, and C.~Liu.
\newblock On the general {L}eslie--{E}ricksen system: Parodi{'}s relation,
  well-posedness and stability.
\newblock {\em Archive for Rational Mechanics and Analysis}, 208:59--107, 2013.

\end{thebibliography}
\end{document}